\documentclass[10pt]{amsart}
\usepackage[pagebackref]{hyperref}

\font\teneur=eurm10
\font\teneus=eusm10

\renewcommand{\d}{\mathrm{d}}
\renewcommand{\L}{\Lambda}
\newcommand{\ts}{\textstyle }
\newcommand{\cF}{{\mathcal F}}
\newcommand{\bbR}{{\mathbb R}}
\newcommand{\bbC}{{\mathbb C}}
\newcommand{\bbP}{{\mathbb P}}
\newcommand{\bbS}{{\mathbb S}}
\newcommand{\bbZ}{{\mathbb Z}}
\newcommand{\ismto}{\,{\tilde{\rightarrow}}\,}
\newcommand{\E}{{\mathrm e}}
\newcommand{\iC}{{\mathrm i}}
\newcommand{\GL}{\operatorname{GL}}
\newcommand{\SL}{\operatorname{SL}}
\newcommand{\SO}{\operatorname{SO}}
\newcommand{\SU}{\operatorname{SU}}
\newcommand{\Or}{\operatorname{O}}
\newcommand{\Spin}{\operatorname{Spin}}
\newcommand{\Scal}{\operatorname{Scal}}
\newcommand{\Ric}{\operatorname{Ric}}
\newcommand{\Ad}{\operatorname{Ad}}

\DeclareMathOperator{\tr}{tr}
\DeclareMathOperator{\End}{End}
\DeclareMathOperator{\Hom}{Hom}
\DeclareMathOperator{\Vect}{Vect}
\DeclareMathOperator{\Diff}{Diff}
\DeclareMathOperator{\Lie}{\hbox{\teneus\char76}}

\newcommand{\phm}{\phantom{-}}
\newcommand{\Gtwo}{\ifmmode{{\rm G}_2}\else{${\rm G}_2$}\fi}
\newcommand{\stars}{\mskip1mu\mathord{*\kern-.5pt_\sigma}\mskip1mu}
\newcommand{\starf}{\mskip1mu\mathord{*\kern-.5pt_\phi}\mskip1mu}

\newcommand{\is}{\mathsf{i}}
\newcommand{\js}{\mathsf{j}}
\newcommand{\Fs}{\mathsf{F}}
\newcommand{\Gs}{\mathsf{G}}

\newcommand{\Qs}{\mathsf{Q}}
\newcommand{\Ss}{\mathsf{S}}
\newcommand{\Vs}{\mathsf{V}}

\newcommand{\euso}{\operatorname{\mathfrak{so}}}
\newcommand{\eugl}{\operatorname{\mathfrak{gl}}}
\newcommand{\eug}{\operatorname{\mathfrak{g}}}
\newcommand{\eut}{\operatorname{\mathfrak{t}}}

\newcommand{\eps}[1]{\varepsilon_{#1}}

\newcommand{\la}{\langle}
\newcommand{\ra}{\rangle}

\newcommand{\ePsi}{\hbox{\teneur\char9}}
\newcommand{\etheta}{\hbox{\teneur\char18}}
\newcommand{\esigma}{\hbox{\teneur\char27}}
\newcommand{\etau}{\hbox{\teneur\char28}}
\newcommand{\epsi}{\hbox{\teneur\char32}}
\newcommand{\eomega}{\hbox{\teneur\char33}}
\newcommand{\eR}{\hbox{\teneur\char82}}
\newcommand{\eS}{\hbox{\teneur\char83}}
\newcommand{\eT}{\hbox{\teneur\char84}}

\newcommand{\w}{{\mathchoice{\,{\scriptstyle\wedge}\,}{{\scriptstyle\wedge}}
      {{\scriptscriptstyle\wedge}}{{\scriptscriptstyle\wedge}}}}
\newcommand{\lhk}{\mathbin{\hbox{\vrule height1.4pt width4pt depth-1pt 
             \vrule height4pt width0.4pt depth-1pt}}}

\numberwithin{equation}{section}

\newtheorem{proposition}{Proposition}
\newtheorem{corollary}{Corollary}

\theoremstyle{remark}
\newtheorem{definition}{Definition}

\newtheorem{remark}{Remark}
\newtheorem{example}{Example}

\begin{document}

\author[R. Bryant]{Robert L. Bryant}
\address{Duke University Mathematics Department\\
         P.O. Box 90320\\
         Durham, NC 27708-0320}
\email{\href{mailto:bryant@math.duke.edu}{bryant@math.duke.edu}}
\urladdr{\href{http://www.math.duke.edu/~bryant}%
         {http://www.math.duke.edu/\lower3pt\hbox{\symbol{'176}}bryant}}

\title[Some remarks on $\Gtwo$-structures]
      {Some remarks on $\Gtwo$-structures}


\begin{abstract}
This article consists of loosely related 
remarks about the geometry of \Gtwo-structures on
$7$-manifolds, some of which are based on unpublished 
joint work with two other people:  F. Reese Harvey and
Steven Altschuler.

After some preliminary background information about
the group \Gtwo\ and its representation theory, a
set of techniques is introduced for calculating the
differential invariants of \Gtwo-structures and
the rest of the article is applications of these
results.  Some of the results that may be of interest 
are as follows: 

First, a formula is derived for the scalar curvature 
and Ricci curvature of a \Gtwo-structure in terms of 
its torsion and covariant derivatives with respect to 
the `natural connection' (as opposed to the Levi-Civita 
connection) associated to a \Gtwo-structure.  When
the fundamental $3$-form of the \Gtwo-structure is 
closed, this formula implies, in particular, that
the scalar curvature of the underlying metric is nonpositive
and vanishes if and only if the structure is torsion-free.
These formulae are also used to generalize a recent
result of Cleyton and Ivanov~\cite{math.DG/0306362} 
about the nonexistence of
closed Einstein $\Gtwo$-structures (other than the 
Ricci-flat ones) on compact $7$-manifolds to a nonexistence
result for closed $\Gtwo$-structures whose Ricci tensor
is too tightly pinched.

Second, some discussion is given of the geometry of
the first and second order invariants of \Gtwo-structures 
in terms of the representation theory of~\Gtwo.

Third, some formulae are derived for closed solutions
of the Laplacian flow that specify how various related
quantities, such as the torsion and the metric, evolve 
with the flow. These may be useful in studying convergence
or long-time existence for given initial data.

Some of this work was subsumed in the work of 
Hitchin~\cite{MR02m:53070} and Joyce~\cite{MR01k:53093}.
I am making it available now mainly because of interest expressed
by others in seeing these results written up since they 
do not seem to have all made it into the literature.
\end{abstract}



\thanks{
Thanks to Rice University and Duke University for their support 
via research grants and to the National Science Foundation 
for its support via DMS-8352009, DMS-8905207, and, most
recently in DMS-0103884.  Thanks also to the
organizers of the April--May 2003 IPAM conference ``Geometry 
and Physics of~$\Gtwo$ manifolds'' for their kind support
and for the interest expressed there in making these
notes available.  Finally, thanks to the organizers
of the 2004 conference on topology and geometry at Gokova, Turkey
for accepting this manuscript to appear in their proceedings
and to their referee for pointing out several typos and mistakes.
\hfill\break
\hspace*{\parindent} 
Version~$1.0$ (\texttt{math.DG/0305124}) 
was posted to the arXiv on 8 May 2003.
This is Version~$5.0$.  
}

\maketitle

\setcounter{tocdepth}{1}
\tableofcontents

\section{Introduction}\label{sec: intro}

This brief article consists of a collection of remarks
on the geometry of $\Gtwo$-structures on $7$-manifolds,
some of which are based on old unpublished joint 
work carried out on separate occasions with two other 
people:  F.~Reese Harvey and Steven Altschuler.  

The work with Reese Harvey (recounted in~\S\ref{sec: torfreecase})
concerned techniques for calculating various quantities 
associated to a \Gtwo-structure, possibly with torsion, 
and was carried out intermittently during the period 1988 
through 1994.  

The work with Steven Altschuler (recounted in~\S\ref{sec: defandevol}) 
concerned the geometry of a natural Laplacian flow for 
\Gtwo-structures and was carried out in 1992.

The main reason for making these remarks available now is 
that some of these formulae and results
do not seem to have appeared yet in the literature and
some people have expressed an interest in learning about 
them.

\section{Algebra}\label{sec: exalg}

This section will collect the main results about the group \Gtwo\
that will be needed.  The reader may consult~\cite{MR89b:53084},
~\cite{MR01k:53093}, or \cite{MR90g:53058} 
for details concerning the properties of the 
group \Gtwo\ that are not proved here.  In general, the notation 
is chosen to agree with the notation in~\cite{MR89b:53084}.

\subsection{The group~\Gtwo}  Let~$e_1,e_2,\ldots,e_7$ denote
the standard basis of~$\bbR^7$ (whose elements will be referred to as 
column vectors of height $7$) and let~$e^1,e^2,\ldots,e^7:\bbR^7\to\bbR$ 
denote the corresponding dual basis.

For notational simplicity, write~$e^{ijk}$ for the wedge 
product~$e^i\w e^j\w e^k$ in~$\L^3\bigl((\bbR^7)^*\bigr)$.  
Define
\begin{equation}\label{eq: phi def}
\phi = e^{123}+e^{145}+e^{167}+e^{246}-e^{257}-e^{347}-e^{356}.
\end{equation}
By a theorem of Schouten~\cite{Schouten1931} 
(see~\cite{MR89b:53084} for a proof), 
the subgroup of~$\GL(7,\bbR)$ that fixes~$\phi$ is a compact, 
connected, simple Lie group of type~\Gtwo. In this article, 
this result will be used to justify the following definition:

\begin{definition}[The group \Gtwo]\label{def: G2}
\begin{equation}\label{eq: G2 def}
\Gtwo = \left\{\ g\in \GL(7,\bbR)\ \vrule\ g^*(\phi) = \phi\ \right\}.
\end{equation}
\end{definition}

\subsection{Associated structures}\label{sssec: assoc strucs}
A few properties of \Gtwo\ will be needed in this article. 
The reader may consult~\cite{MR89b:53084} for proofs.

The group \Gtwo\ acts irreducibly on $\bbR^7$ and preserves the metric 
and orientation for which the basis~$e_1,e_2,\ldots,e_7$ is an oriented 
orthonormal basis.  The notations~$g_\phi$ and~$\la,\ra_\phi$ will be 
used to refer to the metric. The Hodge star operator determined by 
this metric and orientation will be denoted $\starf$.  Note, in 
particular, that \Gtwo\ also fixes the $4$-form
\begin{equation}\label{eq: star phi def}
\starf\phi 
=e^{4567}+e^{2367}+e^{2345}+e^{1357}-e^{1346}-e^{1256}-e^{1247}.
\end{equation}

\subsection{Some \Gtwo\ actions}\label{sssec: Gtwo actions}
The group \Gtwo\ acts transitively on the unit sphere $S^6\subset\bbR^7$.
The stabilizer subgroup of any non-zero vector in~$\bbR^7$ is isomorphic 
to~$\SU(3)\subset\SO(6)$, so that $S^6 = \Gtwo/\SU(3)$.  Since 
$\SU(3)$ acts transitively on~$S^5\subset\bbR^6$, it follows that 
\Gtwo\ acts transitively on the set of orthonormal pairs of vectors 
in~$\bbR^7$. 

However, \Gtwo\ does not act transitively on the set of orthonormal 
triples of vectors in~$\bbR^7$ since it preserves the $3$-form~$\phi$.

\subsection{The $\eps{}$-notation}
It will be convenient to use an $\eps{}$-notation that
will now be introduced.  This is the unique symbol that is 
skew-symmetric in either three or four indices and satisfies 
\begin{align}
 \phi &= {\ts\frac{1}{6}}\eps{ijk}\,e^i\w e^j\w e^k\\
 \starf\phi &= {\ts\frac{1}{24}}\eps{ijkl}\,e^i\w e^j\w e^k\w e^l\,.
\end{align}
Thus, for example, $\eps{123}=1$ and $\eps{4567}=1$, while 
$\eps{124}=\eps{3456}=0$.   Another way to think of this
symbol is via the cross product:  $e_i\times e_j = \eps{ijk}e_k$. 

The symbol $\eps{}$ satisfies various useful identities.  
For example (using the summation convention),
\begin{align}
\eps{ijk}\,\eps{ijl} &=  6\delta_{kl} \label{eq: eps id 1}\\
\eps{ijq}\,\eps{ijkl} &=  4\eps{qkl}  \label{eq: eps id 2}\\
\eps{ipq}\,\eps{ijk} &=
\eps{pqjk}+\delta_{pj}\delta_{qk}-\delta_{pk}\delta_{qj}
     \label{eq: eps id 3}\\
\eps{ipq}\,\eps{ijkl} &=
\delta_{pj}\eps{qkl}-\delta_{jq}\eps{pkl}+
\delta_{pk}\eps{jql}-\delta_{kq}\eps{jpl}+
\delta_{pl}\eps{jkq}-\delta_{lq}\eps{jkp}\,. \label{eq: eps id 4}
\end{align}

These identities are actually quite easy to prove using the fact that 
\Gtwo\ acts transitively on orthonormal pairs.  For example,
identity~\eqref{eq: eps id 3} can be reduced to the case 
where $p=1$ and $q=2$.  
Then the only non-zero term on the left hand side is 
$\eps{312}\eps{3jk}$.  By the definitions of $\phi$ and 
$\starf\phi$, both sides of the equation vanish unless 
$\{j,k\}$ is one of the subsets $\{1,2\}$, $\{4,7\}$, or $\{5,6\}$, and 
the identity clearly holds in those cases.  
The other identities can be proved similarly.

\subsection{Matrix and vector representations}
The $\eps{}$-symbol can be used to describe the algebra $\eug_2$ as 
a subalgebra of $\euso(7)$, the space of skew-symmetric $7$-by-$7$ matrices.
A skew-symmetric matrix $a = (a_{ij})$ lies in~$\eug_2$ if and only if $\eps{ijk}a_{jk}=0$ for all $i$.

For any vector $v=v_{i}e_{i}\in\bbR^7$, define $[v] = 
(v_{ij})\in\euso(7)$  by the formula $v_{ij} = \eps{ijk}v_{k}$. 
It then follows that 
\begin{equation*}
\euso(7) = \eug_2\oplus \left[\bbR^7\right],
\end{equation*}
which is the \Gtwo-invariant irreducible decomposition of 
$\euso(7)$. Note that $[v]$ is the matrix that represents
the linear transformation of~$\bbR^7$ induced by cross-product 
with~$v\in\bbR^7$.

Conversely, define the map $\la\cdot\ra\colon \eugl(7)\to\bbR^7$ 
by $\la\,(a_{ij})\,\ra = (\eps{ijk}a_{jk})$.  The kernel of this 
mapping intersected with~$\euso(7)$ is $\eug_2$ 
and the $\eps{}$-identities imply that, for all 
$a,b\in \bbR^7$,
\begin{align}\label{eq: box and angle identities}
\la\,[a]\,\ra &=6a\\
\la\,[a][b]\,\ra &= 3[b]a = -3[a]b.
\end{align}

\subsection{The \Gtwo-type decomposition of exterior forms}
To avoid writing~$(\bbR^7)^*$ many times, I will, for the rest of
this section, use~$V$ as an abbreviation for the vector space~$\bbR^7$.

Although \Gtwo\ acts irreducibly on $V$ and hence on 
$\L^1(V^*)$ and $\L^6(V^*)$, it does not act 
irreducibly on~$\L^p(V^*)$ for $2\le p\le5$.  In order to 
understand the irreducible decomposition of~$\L^p(V^*)$ for~$p$ 
in this range, it suffices to understand the cases $p=2$ and $p=3$, 
since the operator~$\starf$ induces an isomorphism of 
\Gtwo-modules~$\L^p(V^*)=\L^{7-p}(V^*)$.  

In~\cite{MR89b:53084}, it is shown that there are irreducible 
\Gtwo-module decompositions
\begin{align}
\L^2(V^*) &= \L^2_{14}(V^*)\oplus\L^2_7(V^*) \\
\L^3(V^*) &= 
\L^3_{27}(V^*)\oplus\L^3_7(V^*)\oplus\L^3_1(V^*)
\end{align}
where $\L^p_d(V^*)$ denotes an irreducible \Gtwo-module of 
dimension $d$. For $p$ = $4$ or~$5$, adopt the convention that 
$\L^p_d(V^*)=\starf(\L^{7-p}_d(V^*))$.

These summands can be characterized as follows:
\begin{equation}\label{eq: linear algebra types}
\begin{aligned}
\L^2_7(V^*)
 &=\left\{\,\starf(\alpha\w\starf\phi)\ |\ 
           \alpha\in\L^1(V^*)\ \right\}\\
 &= \left\{\,\alpha\in\L^2(V^*)\ |\ \alpha\w\phi = 
2\starf\alpha\ \right\}\cr
\L^2_{14}(V^*) 
 &= \left\{ \,\alpha\in\L^2(V^*)\ |\ 
     \alpha\w\phi = -\starf\alpha\ \right\} = \eug_2^\flat\\
\noalign{\vskip2pt}
\L^3_1(V^*)
 &= \left\{ \,r\phi\ |\ r\in\bbR\ \right\}\\
\L^3_7(V^*) 
 &=\left\{ \,\starf(\alpha\w\phi)\ |\ \alpha\in\L^1(V^*)\ \right\}\\
\L^3_{27}(V^*) 
 &= \left\{ \,\alpha\in\L^3(V^*)\ |\ 
\alpha\w\phi = 0\ {\rm and}\ \alpha\w\starf\phi = 0\ \right\}
             = \is_\phi\bigl(S^2_0(V^*)\bigr).
\end{aligned}
\end{equation}
The notations $\eug_2^\flat$ and $\is_\phi\bigl(S^2_0(V^*)\bigr)$ 
used in~\eqref{eq: linear algebra types} need some explanation.

First, $\eug_2^\flat$:  Under the ``musical isomorphism'' $\flat\colon 
V\to V^*$ induced by the \Gtwo-invariant inner product 
$\la,\ra_\phi$, the Lie algebra of \Gtwo, namely $\eug_2\subset 
V\otimes V^*$, is identified with $\eug_2^\flat= 
(\flat\otimes1)(\eug_2)\subset \L^2(V^*)\subset V^*\otimes 
V^*$.  This subspace is an irreducible \Gtwo-module since \Gtwo\ is simple.

Second, $\is_\phi\bigl(S^2_0(V^*)\bigr)$: 
Consider the linear mapping $\is_\phi\colon S^2(V^*) \to \L^3(V^*)$, 
defined on decomposable elements by 
\begin{equation}\label{eq: is def}
\is_\phi(\alpha\circ\beta) = 
\alpha\w\starf(\beta\w\starf\phi)
+\beta\w\starf(\alpha\w\starf\phi).
\end{equation}
The mapping~$\is_\phi$ is \Gtwo-invariant and one can show 
that~$S^2(V^*) = \bbR g_\phi \oplus S^2_0(V^*)$ is a decomposition 
of~$S^2(V^*)$ into \Gtwo-irreducible summands.  Evidently,
$\is_\phi$ is nonzero on each summand and is therefore injective.  
Hence, the image~$\is_\phi\left(S^2_0(V^*)\right)\subset\L^3(V^*)$ 
is $27$-dimensional and irreducible.  The equation
\begin{equation}
\L^3_{27}(V^*) 
 = \left\{ \,\alpha\in\L^3(V^*)\ |\ 
\alpha\w\phi = 0\ {\rm and}\ \alpha\w\starf\phi = 0\ \right\}
\end{equation}
defines~$\L^3_{27}(V^*)$ as a \Gtwo-invariant, $27$-dimensional 
subspace of~$\L^3(V^*)$.  By dimension count, it must
intersect~$\is_\phi\left(S^2_0(V^*)\right)$ nontrivially.
Since this intersection is also \Gtwo-invariant and
since~$\is_\phi\left(S^2_0(V^*)\right)$ is 
\Gtwo-irreducible,~$\is_\phi\left(S^2_0(V^*)\right)=\L^3_{27}(V^*)$.

Using the $\varepsilon$-notation, one can express the map~$\is_\phi$ 
in indices as
\begin{equation}\label{eq: is in indices}
\is_\phi(h_{ij}e^ie^j) = \varepsilon_{ikl}\,h_{ij}\,e^j\w e^k\w e^l,
\end{equation}
making it evident that~$\is_\phi(g_\phi) = 6\phi$.  

It will be useful to have a way to invert the map~$\is_\phi$.  
Define~$\js_\phi:\L^3(V^*)\to S^2(V^*)$
by the formula
\begin{equation}\label{eq: js def}
\js_\phi(\gamma)(v,w) = \starf\bigl((v\lhk\phi)\w(w\lhk\phi)\w\gamma\bigr).
\end{equation}
for~$\gamma\in \L^3(V^*)$ and~$v,w\in V$.  It is not difficult
to verify that
\begin{equation}\label{eq: js comp is}
\js_\phi\bigl(\is_\phi(h)\bigr) = 8h + 4\bigl(\tr_{g_\phi}(h)\bigr)\,g_\phi
\end{equation}
for all~$h\in S^2(V^*)$.  Note also that $\js_\phi(\phi) = 6 g_\phi$, 
while~$\js_\phi\bigl(\L^3_7(V^*)\bigr) = 0$.

Note that~$\is_\phi$ and~$\js_\phi$ are not isometries
when~$S^2_0(V^*)$ and~$\L^3_{27}(V^*)$ are given their natural 
metrics.\footnote{The usual inner product on
exterior forms is meant here, while, when~$h = h_{ij}\,e^ie^j$
with~$(e^i)$ being a $g$-orthonormal coframe of~$V$, one sets
$|h|^2 = h_{ij}h_{ij}$.} 
Instead,~$\gamma\in\L^3_{27}(V^*)$ satifies~$|\js_\phi(\gamma)|^2 
= 8\,|\gamma|^2$ while ~$h\in S^2_0(V^*)$ 
satisfies~$|\is_\phi(h)|^2 = 8\,|h|^2$. 

\subsection{More \Gtwo\ representation theory}
It will, from time to time, be useful to have some deeper
knowledge of the representation theory of~$\Gtwo$, so some
of these facts will be collected here.  For details, 
consult~\cite{MR81b:17007}.

Since~$\Gtwo$ is a simple Lie group of rank~$2$, its irreducible
representations can be indexed by a pair of integers~$(p,q)$ 
that represent the highest weight of the representation with
respect to a fixed maximal torus in~$\Gtwo$ endowed with fixed
base for its root system.  The irreducible representation 
of highest weight~$(p,q)$ will be denoted~$\Vs_{p,q}$. 

\subsubsection{The standard representation}\label{sssec: R7 as G2rep}
The fundamental representation~$\Vs_{1,0}\simeq\bbR^7$ is the 
`standard' representation in which~\Gtwo\ has been defined
in this article.

The representation~$\Vs_{p,0}$ for~$p\ge0$ is isomorphic
to $\Ss^p_0(\bbR^7)$, i.e., the symmetric, trace-free polynomials
of degree~$p$ in seven variables.  (It is somewhat remarkable
that these irreducible representations of~$\SO(7)$ remain
irreducible when thought of as representations of~$\Gtwo$.)
In this article, the only representations~$\Vs_{p,0}$ in
this series that will be important are those for~$p=0,1,2$.

\subsubsection{The adjoint representation}\label{sssec: adjoint repn}
The other fundamental representation,~$\Vs_{0,1}\simeq\bbR^{14}$
is isomorphic to~$\eug_2$, i.e., is the adjoint representation
of~$\Gtwo$.  The representation~$\Vs_{0,p}$ for~$p\ge0$ is then the
irreducible constituent of highest weight in~$\Ss^p(\eug_2)$.

In this article, only~$\Vs_{0,1}\simeq\eug_2$ and~$\Vs_{0,2}
\simeq\bbR^{77}$ from this series will be important.  (This
latter one will be important because it is the space of 
curvature tensors of $\Gtwo$-metrics.)  The reader must be
careful not to confuse the representation~$\Vs_{0,2}$ with
$\Vs_{3,0}$, which also happens to have dimension~$77$.

A few more facts about this representation
will be needed:  The group~$\Gtwo$ has rank~$2$ and
a maximal torus for~$\Gtwo$ can be obtained by simply
taking a maximal torus in the subgroup~$\SU(3)$.  
Moreover, every element in~$\eug_2$ is $\Ad(\Gtwo)$-conjugate
to an element in such a maximal torus.  Consequently,
every element in~$\L^2_{14}(\bbR^7)=\eug_2^\flat$ is
conjugate to an element of the form
\begin{equation}\label{eq: toral form for g2}
\alpha = \lambda_1\,e^{23} + \lambda_2\,e^{45}
          - (\lambda_1{+}\lambda_2)\,e^{67}
\end{equation}
since these span~$\eut^\flat\subset\eug_2^\flat$,
where~$\eut\subset\eug_2$ is a Cartan subalgebra.  Moreover, 
it is well-known that the ring of $\Ad(\Gtwo)$-invariant polynomials 
on~$\eug_2$ is a free polynomial ring on two generators, 
one of degree~$2$ and one of degree~$6$.  One sees from
the above normal form that these two generators can be taken
to be $|\alpha|^2$ and $|\alpha^3|^2$.  Thus, two elements
$\alpha$ and~$\beta$ in~$\L^2_{14}(\bbR^7)$ are 
conjugate under the action of~$\Gtwo$ if and only if 
they satisfy~$|\alpha|^2=|\beta|^2$ and~$|\alpha^3|^2=|\beta^3|^2$.
In particular, the normal form~\eqref{eq: toral form for g2} 
can be made unique by requiring that~$0\le\lambda_1\le\lambda_2$.

In particular, one obtains, 
for all $\alpha\in\L^2_{14}(V^*)$, the useful identity
\begin{equation}\label{eq: G2 quartic}
|\,\alpha^2\,|^2 = |\alpha|^4
\end{equation}
and inequality
\begin{equation}\label{eq: G2 sextic ineq}
|\,\alpha^3\,|^2 \le {\ts\frac23}\,|\alpha|^6,
\end{equation} 
which are easily verified by checking them on 
elements of the form~\eqref{eq: toral form for g2}. 

In fact, using the normal form~\eqref{eq: toral form for g2}, 
one can prove other useful exterior algebra identities.  
One that will be needed later is
\begin{equation}\label{eq: G2 cubic express}
\alpha\w\starf(\alpha\w\alpha) = |\alpha|^2\,\starf\alpha
-{\ts\frac13}\,\starf(\alpha^3)\w\starf\phi
\qquad \text{for $\alpha\in\L^2_{14}(V^*)$.}
\end{equation}

\subsubsection{Other representations}\label{sssec: other repns}
Of the representations~$\Vs_{p,q}$ with~$p$ and~$q$ positive,
only~$\Vs_{1,1}\simeq\bbR^{64}$ will play any significant role
in this article (and mainly as a nuisance at that).  In fact,
each of the other representations~$\Vs_{p,q}$ with both $p$ and~$q$ 
positive has dimension at least~$189$, so these can easily be 
ruled out for dimension reasons in the calculations to follow.  

The following tensor product and Schur functor decompositions
will be useful: 
\begin{equation}\label{eq: tensor and schur decomps}
\begin{aligned}
\Ss^2\bigl(V_{1,0}\bigr) &\simeq\Vs_{0,0}\oplus\Vs_{2,0}\\
\L^2\bigl(V_{1,0}\bigr) &\simeq\Vs_{1,0}\oplus\Vs_{0,1}\\
\Vs_{1,0}\otimes \Vs_{0,1} &\simeq\Vs_{1,0}\oplus\Vs_{2,0}\oplus\Vs_{1,1}\\
\Ss^2\bigl(V_{0,1}\bigr) &\simeq\Vs_{0,0}\oplus\Vs_{2,0}\oplus\Vs_{0,2}\\
\L^2\bigl(V_{0,1}\bigr) &\simeq\Vs_{0,1}\oplus\Vs_{3,0}
\end{aligned}
\end{equation}

\subsubsection{An example of $\Gtwo$-type decomposition}
\label{sssec: G2 decomp exmpl}
As an application of these formulae that will be used below,
consider the problem of decomposing~$\beta\w\beta\in\L^4(V^*)$ 
into its $\Gtwo$-types where~$\beta$ lies in~$\L^2_{14}(V^*)\simeq
\Vs_{0,1}$. Since
\begin{equation}
\L^4(V^*)\simeq
\L^4_1(V^*)\oplus\L^4_7(V^*)\oplus\L^4_{27}(V^*)
\simeq\Vs_{0,0}\oplus\Vs_{1,0}\oplus\Vs_{0,2}
\end{equation}
and since, by~\eqref{eq: tensor and schur decomps},
we have~$\Ss^2\bigl(\Vs_{0,1}\bigr)
 \simeq\Vs_{0,0}\oplus\Vs_{2,0}\oplus\Vs_{0,2}$,
it follows that $\beta\w\beta$ can have no component in $\L^4_7(V^*)
\simeq \Vs_{1,0}$.  Moreover, since there is, up to multiples,
only one $\Gtwo$-invariant quadratic form on~$\Vs_{0,1}$ and since
$\starf\phi$ spans~$\L^4_1(V^*)\simeq \Vs_{0,0}$, it follows that
there is a constant~$\lambda$ such that
\begin{equation}
\beta\w\beta = \lambda\,|\beta|^2\,\starf\phi 
+ \bigl(\beta\w\beta - \lambda\,|\beta|^2\,\starf\phi\bigr)
\end{equation}
where the first term on the right lies in~$\L^4_1(V^*)$ while
the second term (in parentheses) lies in~$\L^4_{27}(V^*)$.

The constant~$\lambda$ is determined as follows:  Wedging both
sides with~$\phi$ and using the fact that~$\beta\w\phi = -\starf\beta$
while $\gamma\w\phi=0$ for~$\gamma\in\L^4_{27}(V^*)$
yields
\begin{equation}
-|\beta|^2\starf1 = \beta\w\beta\w\phi 
= \bigl(\lambda\,|\beta|^2\,\starf\phi\bigr)\w\phi 
= 7\lambda\,|\beta|^2\starf1,
\end{equation}
showing that~$\lambda = -\frac17$.  
Thus, the $\Gtwo$-type decomposition is given by
\begin{equation}\label{eq: beta2 decomp}
\beta\w\beta = -{\ts\frac17}\,|\beta|^2\,\starf\phi 
+ \bigl(\beta\w\beta + {\ts\frac17}\,|\beta|^2\,\starf\phi\bigr).
\end{equation}
for~$\beta\in\L^2_{14}(V^*)$.  Of course, this decomposition
is orthogonal, so, using the identity~\eqref{eq: G2 quartic}, 
one can take the square norms of both sides, yielding
\begin{equation}
|\beta|^4 = |\beta\w\beta|^2 = {\ts\frac17}|\beta|^4
+ \bigl|\beta\w\beta + {\ts\frac17}\,|\beta|^2\,\starf\phi\bigr|^2.
\end{equation}
Consequently, for~$\beta\in\L^2_{14}(V^*)$, one has
\begin{equation}\label{eq: 27-piece quartic id}
\bigl|\,\beta\w\beta + {\ts\frac17}\,|\beta|^2\,\starf\phi\,\bigr|^2
= {\ts\frac67}|\beta|^4, 
\end{equation}
an identity that will be used below.  
(Note that~\eqref{eq: 27-piece quartic id} implies,
in particular, that the $\L^4_{27}(V^*)$-piece of~$\beta\w\beta$
cannot vanish unless~$\beta$ itself vanishes, a result equivalent to
Lemma~5.8 of~\cite{math.DG/0306362}.)

Similar sorts of calculations can be used to establish the
(sharp) inequalities for quadratic forms
\begin{equation}\label{eq: jbetabeta bounds}
-2\,|\beta|^2\,g\ \le\ \js\bigl(\starf(\beta\w\beta)\bigr) 
\ \le\ {\ts\frac23}\,|\beta|^2\,g.
\end{equation}
Details are left to the reader.

\subsection{Definite forms}\label{sssec: definite forms}
The dimension of \Gtwo\ is~$14$ and so, by dimension count,  
the $\GL(V)$-orbit of~$\phi$ in~$\L^3(V^*)$ is open.  
Denote this orbit by $\L^3_+(V^*)$ and speak of the 
elements of $\L^3_+(V^*)$ as \emph{definite} $3$-forms 
on~$V$.  Note that $\L^3_+(V)$ has two components, since~$\GL(V)$
does and since~\Gtwo\ is connected.  Each component is the 
negative of the other.  It is known~\cite{MR91e:53056} that
$\SO(7)/\Gtwo\simeq\bbR\bbP^7$, so that each component 
of~$\L^3_+(V)$ is diffeomorphic to~$\bbR\bbP^7\times\bbR^{28}$.

\subsubsection{On general $7$-dimensional vector spaces}
If $W$ is any $7$-dimensional vector space, an 
isomorphism $u\colon W\ismto V$ induces an isomorphism 
$u^*\colon\L^3(V^*)\ismto\L^3(W^*)$. 
Denote by~$\L^3_+(W^*)$ the open subset 
$u^*\left(\L^3_+(V^*)\right)\subset\L^3(W^*)$.  
Since~$\L^3_+(V^*)$ consists of a single~$\GL(V)$-orbit, 
this set does not depend on the choice of~$u$.

\subsubsection{Associated algebraic structures}
Each~$\varphi\in\L^3_+(W^*)$ has a stabilizer in~$\GL(W)$
that is isomorphic to~$\Gtwo$ and hence defines a canonical 
inner product~$\la,\ra_\varphi$ (with associated quadratic
form~$g_\varphi$) and orientation 
(Hodge star)~$\ast_\varphi:\L^p(W^*)\to\L^{7-p}(W^*)$.

Similarly, using~$\varphi$ in the place of~$\phi$ in the
formulae~\eqref{eq: is def} and~\eqref{eq: js def}, one
defines mappings~$\is_\varphi:S^2(W^*)\to\L^3(W^*)$
and~$\js_\varphi:\L^3(W^*)\to S^2(W^*)$.  These maps are
frequently useful in formulae.  

For example, let~$\Gs:\L^3_+(W^*)\to S^2_+(W^*)$ be the
nonlinear $\GL(W)$-equivariant mapping that 
satisfies~$\Gs(\varphi) = g_\varphi$.  It is not difficult
to show that~$\Gs$ is smooth and satisfies
\begin{equation}
\Gs'(\varphi)(\psi) = {\ts\frac12}\,\js_\varphi(\psi)
-{\ts\frac13}\,{\ast_\varphi}(\psi\w{\ast_\varphi}\varphi)\,g_\varphi\,.
\end{equation}

There is also an associated \emph{vector cross product}
$\times_\varphi: W\times W \to W$ defined by the condition
\begin{equation}
\la w_1\times_\varphi w_2\,, w_3\ra_\varphi
 = \varphi(w_1,w_2,w_3).
\end{equation}

\begin{remark}[The vector cross product definition of \Gtwo]
Given a vector space~$V$ over~$\bbR$ endowed with 
a positive definite inner product $\la,\ra:V\times V\to\bbR$, 
a ($2$-fold) \emph{vector cross product} on~$\bigl(V,\la,\ra\bigr)$
is a skew-symmetric bilinear pairing~$\times:V\times V\to V$
that satisfies
\begin{equation}
\la v_1\times v_2, v_1\ra = 0
\qquad\text{and}\qquad
|v_1\times v_2|^2 = |v_1|^2\,|v_2|^2 - \la v_1,v_2\ra^2
\end{equation}
for all~$v_1,v_2\in V$.  

It can be shown that the $\GL(7,\bbR)$-stabilizer of the vector
cross product~$\times_\phi$ is equal to~\Gtwo.  
Hence one could take~$\times_\phi:\bbR^7\times\bbR^7\to\bbR^7$ 
as the algebraic structure defining \Gtwo.  In fact, this
is what Gray did in his work on \Gtwo-structures.  However,
I find that the $3$-form formulation is more congenial for
computations, so vector cross products will not play any significant
role in this article.
\end{remark}

\subsubsection{Definite $4$-forms}
The canonical mapping~$\Ss:\L^3_+(W^*)\to\L^4(W^*)$
defined by~$\Ss(\varphi) = \ast_\varphi\varphi$ is 
a double covering onto an open set~$\L^4_+(W^*)$
in~$\L^4(W^*)$, which will be referred to as the
space of `definite' $4$-forms on~$W$.  

The $\GL(W)$-stabilizer of an element~$\psi\in\L^4_+(W^*)$ 
is then isomorphic to~$\pm\Gtwo 
= \Gtwo\cup \bigl(\Gtwo\cdot(-\text{id}_W)\bigr)$.
Thus, a definite $4$-form on~$W$ defines an inner product 
on~$W$, but not an orientation.

\section{\Gtwo-structures}  

\subsection{Definite forms on manifolds}
Let $M$ be a smooth manifold of dimension 7. 
The union of the subspaces $\L^3_+(T_x^*M)$ 
is an open subbundle $\L^3_+(T^*M)\subset 
\L^3(T^*M)$ of the bundle of $3$-forms on~$M$.  

\begin{definition}[Definite $3$-forms on manifolds]
A $3$-form~$\sigma$ on~$M$ that takes values in~$\L^3_+(T^*M)$
will be said to be a \emph{definite} $3$-form on~$M$.
The set of definite $3$-forms on~$M$ will be denoted $\Omega^3_+(M)$.
\end{definition}

\subsubsection{\Gtwo-structures and definite $3$-forms}
Each definite~$3$-form on~$M$ defines a \Gtwo-structure on~$M$ 
in the following way:

Let $\cF$ denote the principal right $\GL(V)$-bundle over $M$ 
consisting of $V$-coframes $u\colon T_xM\ismto V$. Given any 
$\sigma\in\Omega^3_+(M)$, define a \Gtwo-bundle
\begin{equation}\label{eq: Fsigma def}
F_\sigma
=\left\{u\in\Hom(T_xM,V)\ \mid\ x\in M\ {\rm and}\ u^*(\phi)=\sigma_x\right\}.
\end{equation}

Every \Gtwo-reduction of~$\cF$ (i.e., \Gtwo-structure on~$M$
in the usual sense) is of the form~$F_\sigma$ for some unique~$\sigma\in\Omega^3_+(M)$.  For this reason, 
a $3$-form~$\sigma\in\Omega^3_+(M)$ will usually, 
by abuse of language, be called a \Gtwo-structure in this article.

\begin{remark}[Alternative terminologies]\label{rem: altterms}
Some authors use `almost \Gtwo-structure' to refer
to what I am calling a \Gtwo-structure in this article. 
Apparently, this practice stems from an imagined analogy with the
distinction between `almost complex structure' and `complex structure'.

However, for a subgroup~$G\subset\GL(n,\bbR)$, the use 
of `$G$-structure' on an $n$-manifold~$M$ to mean a $G$-subbundle 
of the $\GL(n,\bbR)$-bundle of frames (or coframes) on~$M$
is well established.  It seems unwise to tamper with this
usage, especially since `almost $G$-structure' 
suggests a structure that lacks some property of
actual $G$-structures.  Making an exception 
for the case~$G=\Gtwo$ merely invites confusion.

This use of `$G$-structure' does not conflict with the 
`almost complex structure' vs.~`complex structure' usage 
since a complex structure on a $2n$-manifold is not simply 
a $\GL(n,\bbC)$-structure, but is (by the Newlander-Nirenberg 
theorem, equivalent to) a $\GL(n,\bbC)$-structure 
with an assumed integrability property,
whereas an `almost complex structure' actually is 
(equivalent to) a $\GL(n,\bbC)$-structure, 
not an `almost $\GL(n,\bbC)$-structure'.

Some authors speak of an `integrable \Gtwo-structure', meaning
a \Gtwo-structure~$\sigma\in\Omega^3_+(M)$ satisfying some 
differential equations, such as~$\d\sigma=0$ (the exact 
differential equation intended varies with the author).  
Again, this usage appears to stem from an imagined analogy with
a symplectic structure, which is defined by a nondegenerate 
$2$-form~$\omega$ that is closed, i.e., $\d\omega=0$.
In the symplectic case, Darboux' Theorem says that~$\omega$ 
is, indeed, locally equivalent to the flat model, i.e., is 
`integrable' in the standard terminology of the theory of 
Lie pseudo-groups.  (In a similar way, one speaks of `integrable 
almost complex structures'.)  This usage of `integrable' 
for \Gtwo-structures also seems ill-advised 
to me since, as will be seen below, no first order condition
on a \Gtwo-structure implies that it is locally equivalent 
to the flat model (which is the only interpretation of 
`integrable' in this context that would be consistent with
the established usage in the theory of Lie pseudo-groups).  
Moreover, this encourages the confusing shift of terminology
in which `\Gtwo-structure' is used to mean `integrable \Gtwo-structure'
and `almost \Gtwo-structure' is used to mean an actual 
\Gtwo-structure.  

For this reason, none of the modifiers `integrable', `almost',
`nearly', or their ilk will be used in this article when 
referring to \Gtwo-structures.  

However, since it seems to be harmless, the terminology 
`\Gtwo-manifold' will sometimes be used to denote a manifold
endowed with a \Gtwo-structure that is flat to first order
(i.e., `torsion-free' in the usual terminology).  
\end{remark}

\begin{definition}[Associated metric, orientation, 
and vector cross product]
\label{def: gsigma def}
For any~$\sigma\in\Omega^3_+(M)$, denote by $g_\sigma$, $*_\sigma$,
and~$\times_\sigma$ the metric, Hodge star operator, and
vector cross product on~$M$ that are canonically associated to~$\sigma$.  
When it is needed, the oriented orthonormal frame bundle of~$g_\sigma$
with this orientation will be denoted~$\Fs_\sigma = F_\sigma\cdot\SO(7)$.
\end{definition}

\begin{remark}[Existence of \Gtwo-structures]
\label{rem: G2existence}
Because \Gtwo\ is both connected and simply connected,
a connected $7$-manifold~$M$ can support a \Gtwo-structure 
only if it is both orientable and spinnable, i.e., if the first
two Stiefel-Whitney classes of~$M$ vanish.

Conversely, by an observation due to Gray~\cite{MR39:4790}, 
these two necessary conditions are also sufficient:  

Since~$\Gtwo$ is simply connected, it is the image under 
the standard double covering map~$\rho:\Spin(7)\to\SO(7)$ 
of a unique subgroup of $\Spin(7)$, which, by abuse of language, 
will also be called~$\Gtwo$.  Now, $\Spin(7)$ has a faithful 
representation on~$\bbR^8$ and hence can be regarded as
a subgroup of~$\SO(8)$.  The restriction of this representation
to~$\Gtwo$ must also be faithful and hence, for dimension
reasons, it must be isomorphic to $\Vs_{0,0}\oplus\Vs_{1,0}$.
In particular, $\Gtwo$ fixes a vector in~$\bbR^8$ and
acts transitively on the unit $6$-sphere orthogonal to
this vector.  Consequently, $\Spin(7)$ must act transitively
on the unit $7$-sphere in~$\bbR^8$ with stabilizer 
subgroup~$\Gtwo$.

Now, suppose~$M^7$ to be orientable and spinnable.  
Choose a Riemannian metric~$g$, an orientation,
and a spin structure~$\tilde \Fs\to M$,
i.e., a spin double cover of the $\SO(7)$-bundle~$\Fs\to M$ 
consisting of oriented, $g$-orthonormal coframes on~$M$.
The associated spinor bundle~$\bbS=\tilde\Fs\times_{\Spin(7)}\bbR^8$
is a vector bundle of rank~$8$ over the $7$-manifold~$M$
and therefore has a nonvanishing unit section~$s:M\to\bbS$. 
This allows one to reduce the structure group of~$\tilde\Fs$
(and hence~$\Fs$) from~$\Spin(7)$ to~$\Gtwo$ (since, by the 
previous paragraph, this is, up to conjugacy, the 
$\Spin(7)$-stablizer of any nonzero vector in~$\bbR^8$). 
Thus, $M$ admits a $\Gtwo$-structure whose associated
metric and orientation are the chosen ones. 
\end{remark}

\subsection{Type decomposition}
Since \Gtwo\ acts reducibly on $\L^p(V^*)$ for $2\le p\le 5$, 
one can associate to any \Gtwo-structure $\sigma$ on~$M$ natural 
splittings of the $p$-form bundles $\L^p(T^*M)$ into 
direct summands.  These will be labeled as $\L^p_d(T^*M,\sigma)$, 
or more simply, $\L^p_d(T^*M)$ when the structure $\sigma$ is 
clear from context.  Denote the space of sections of 
$\L^p_d(T^*M,\sigma)$ by $\Omega^p_d(M,\sigma)$.  

Thus, for example, in view of~\eqref{eq: linear algebra types}, 
one has  
\begin{align}\label{eq: 2-form types}
\Omega^2_7(M,\sigma) &= \left\{\, \beta\in\Omega^2(M)\ |\ 
\beta\w\sigma=2\stars\beta\ \right\}\\
\Omega^2_{14}(M,\sigma) &= \left\{\, \beta\in\Omega^2(M)\ |\ 
\beta\w\sigma=-\stars\beta\ \right\}.
\end{align}

Fortunately, the irreducible modules of dimensions~$14$ and~$27$
only occur in one dual pair of dimensions each.  Meanwhile, the
irreducible module of dimension~$7$ occurs in each 
degree~$1\le p\le 6$.  From time to time, it is useful to be able
to recognize the scale factors that can be introduced by the
various different isomorphisms between these different modules.
For example, for~$\alpha\in\Omega^1_7(M)$ one has
\begin{equation}\label{eq: star ids for 1-forms}
\begin{aligned}
\stars\bigl(\stars(\alpha\w\sigma)\w\sigma\bigr) &= -4\,\alpha\\
\stars\bigl(\stars(\alpha\w\stars\sigma)\w\stars\sigma\bigr) &= \phm3\,\alpha,
\end{aligned}
\end{equation}
and these identities can sometimes be useful in simplifying 
various expressions.  One should also keep in mind that, 
using the metric, each $1$-form~$\alpha$ has a corresponding
dual vector field~$\alpha^\sharp$ and there are useful
identities of the form
\begin{equation}\label{eq: vectors and $1$-forms}
\begin{aligned}
\stars(\alpha\w\sigma) &= -\alpha^\sharp \lhk \stars\sigma\\
\stars(\alpha\w\stars\sigma) &= \phm\alpha^\sharp \lhk \sigma\,.
\end{aligned}
\end{equation}

\begin{remark}[\Gtwo-structures with the same associated metric
and orientation]\label{rem: sameassmetandorient}
These type decompositions have many uses.  For example, 
they furnish a description of all of the
\Gtwo-structures that have the same associated metric and
orientation as a given~$\sigma\in\Omega^3_+(M)$:

Let~$a$ and~$\alpha$ be a function and a $1$-form, respectively,
on~$M$ with~$a^2 + |\alpha|^2_\sigma = 1$.  Then the $3$-form
\begin{equation}\label{eq: tsigmainaalpha}
\tilde \sigma 
= \bigl(a^2 - |\alpha|^2_\sigma)\,\sigma
   + 2 a\,\stars\bigl(\alpha \w \sigma)
   + \is\bigl(\alpha\circ\alpha)
\end{equation}
is definite and has the same associated metric and orientation
as~$\sigma$. (This pointwise fact is most easily proved by
checking it in the case~$\sigma=\phi$ and~$(a,\alpha) = (c,s\,e^1)$
where~$c^2+s^2=1$ and then using the fact that~\Gtwo\ acts transitively
on the unit $6$-sphere in~$\bbR^7$ to reduce to this case.) 

Moreover, any definite $3$-form on~$M$ that
has~$g_\sigma$ and~$\stars$ as associated metric and orientation
is of the form~\eqref{eq: tsigmainaalpha} 
for some pair~$(a,\alpha)$ satisfying $a^2 + |\alpha|^2_\sigma = 1$, 
unique up to replacement by~$(-a,-\alpha)$.  
(If~$H^1(M,\bbZ_2)\not=0$, 
the pair~$(a,\alpha)$ might only be defined up to sign.)

Of course, some such formula was expected, 
since~$\SO(7)/\Gtwo\simeq \bbR\bbP^7$ (a consequence
of the result~$\Spin(7)/\Gtwo\simeq S^7$ discussed 
in Remark~\ref{rem: G2existence}).  What~\eqref{eq: tsigmainaalpha} 
displays is a concrete isomorphism between the bundle~$\Fs_\sigma/\Gtwo$
and the~$\bbR\bbP^7$-bundle~$\bbP\bigl(\bbR\oplus T^*M\bigr)$ over~$M$.  
\end{remark}

\subsection{Exterior derivative formulae}
The decomposition of the $p$-forms on $M$ allows one to express the 
exterior derivatives of both $\sigma$ and $\stars\sigma$ in fairly 
simple terms:

\begin{proposition}[The torsion forms]\label{prop: ext der formulae}
For any \Gtwo-structure $\sigma\in\Omega^3_+(M)$, there 
exist unique differential forms $\tau_0\in\Omega^0(M)$, 
$\tau_1\in\Omega^1(M)$, $\tau_2\in\Omega^2_{14}(M,\sigma)$, and 
$\tau_3\in\Omega^3_{27}(M,\sigma)$ so that the following equations 
hold:
\begin{equation}\label{eq: dsigma delsigma}
\begin{aligned}
\d\,\sigma&=\tau_0\,\stars\sigma+3\,\tau_1\w\sigma+\stars\tau_3\,,\\
\d\,\stars\sigma &=\qquad4\,\tau_1\w\stars\sigma+ 
\tau_2\w\sigma\,.
\end{aligned}
\end{equation}
\end{proposition}

\begin{proof} In view of the 
decomposition~\eqref{eq: linear algebra types},
the only part of this proposition that is not 
simply the definition of the $\tau_i$ is the occurrence of~$\tau_1$ 
in two places.  In fact, by~\eqref{eq: linear algebra types}, 
there exist unique forms 
$\tau_0\in\Omega^0(M)$, $\tau_1,\tilde\tau_1\in\Omega^1(M)$, 
$\tau_2\in\Omega^2_{14}(M,\sigma)$, and 
$\tau_3\in\Omega^3_{27}(M,\sigma)$
so that the above equation for $d\,\sigma$ holds while
$\d\,\stars\sigma=4\,\tilde\tau_1\w\stars\sigma+\tau_2\w\sigma$.  

However, as is shown in~\cite{MR89b:53084} 
(see Remark~\ref{rem: geninttor} below 
for a sketch of the proof), there is an identity
\begin{equation}\label{eq: torsion identity}
\stars\sigma\w\stars\bigl(\d(\stars\sigma)\bigr) + 
(\stars\d\sigma)\w\sigma = 0
\end{equation}
valid for all $\sigma\in\Omega^3_+(M)$, 
and, in view of~\eqref{eq: star ids for 1-forms},
this is equivalent to $\tilde\tau_1 = \tau_1$. 
\end{proof}

\begin{definition}[The torsion forms]\label{def: torsion forms}
For a definite~$3$-form~$\sigma\in\Omega^3_+(M)$, the quadruple
of forms~$(\tau_0,\tau_1,\tau_2,\tau_3)$ defined 
by~\eqref{eq: dsigma delsigma} will be referred to as
the \emph{intrinsic torsion forms} of~$\sigma$.
\end{definition}

\begin{remark}[General intrinsic torsion]\label{rem: geninttor}
The existence of the identity~\eqref{eq: torsion identity}
may seem surprising at first, but the existence of such an identity
can be understood by general considerations. 

For any subgroup~$G\subset\SO(n)$, the first order invariants 
(usually called the `intrinsic torsion') of a~$G$-structure~$F$ 
on an $n$-manifold~$M$ take values in a bundle over~$M$ associated 
to the natural~$G$-representation on~$(\euso(n)/\eug)\otimes\bbR^n$.
(See \S\ref{ssec: bundle torsion} below for a further explication
of this fact.)  When the first order invariants of a given
$G$-structure vanish, it is said to be `$1$-flat' 
or `flat to first order'.  For more discussion of this notion,
see~\cite{MR89b:53084}.
\end{remark}

In the case of~$\Gtwo\subset\SO(7)$, 
this torsion representation space is
\begin{equation}\label{eq: g2 torsion}
(\euso(7)/\eug_2)\otimes\bbR^7\simeq \Vs_{1,0}\otimes\Vs_{1,0}
\simeq \Vs_{0,0}\oplus\Vs_{1,0}\oplus\Vs_{0,1}\oplus\Vs_{2,0}. 
\end{equation}
and, as has already been remarked, these four summands are isomorphic, 
respectively, to~$\L^0(V^*)$, $\L^1(V^*)$, $\L^2_{14}(V^*)$,
and~$\L^3_{27}(V^*)$.  Since the exterior derivatives of the
defining forms~$\sigma$ and~$\stars\sigma$ can be expressed linearly
in terms of the first order invariants of~$F_\sigma$ and since
there is only one~$\L^1(V^*)$ in the above representation list, 
it follows that the two $1$-forms~$\tau_1$ and~$\tilde\tau_1$ 
alluded to in the above proof must satisfy some universal linear relation.

Consideration of the fact that replacing $\sigma$ by~$\lambda^3\sigma$
for some positive function~$\lambda$ will replace~$\stars\sigma$
by~$\lambda^4\stars\sigma$ shows that this relation must be
the one given in Proposition~\ref{prop: ext der formulae}. 

\begin{proposition}[$1$-flatness of \Gtwo-structures]
\label{prop: 1flatg2}
A \Gtwo-structure~$\sigma\in\Omega^3_+(M)$ is flat to first
order if and only if its torsion forms all vanish, i.e., if
and only if $\d\sigma = \d\stars\sigma=0$.
\end{proposition} 

\begin{proof}
A \Gtwo-structure~$\sigma\in\Omega^3_+(M)$ is flat to first 
order at~$p\in M$ if there exists a $p$-centered coordinate
chart~$x:U\to\bbR^7$ such that the $3$-form~$\sigma - x^*(\phi)$
on~$U$ vanishes to order at least~$2$ at~$p$.  

Recall that the map~$\Ss:\L^3_+(W^*)\to\L^4_+(W^*)$ 
defined in \S\ref{sssec: definite forms} is a smooth double 
covering.  This implies that if~$\sigma - x^*(\phi)$ vanishes 
to order~$2$ at~$p$, then $\stars\sigma - x^*(\ast_\phi\phi)$ 
vanishes to order~$2$ at~$p$ as well.  

Since~$\d\phi = \d{\ast_\phi}\phi = 0$, if~$\sigma$ is flat 
to first order at~$p$, then $\d\sigma$ and~$\d\stars\sigma$ 
must vanish to at least first order at~$p$.
Thus, the claim in one direction is established. 

To demonstrate the claim in the converse direction, 
it suffices to show that any definite $3$-form~$\sigma$ 
defined on a neighborhood of~$0\in\bbR^7$ 
that satisfies~$\sigma_0 = \phi$ and~$\d\sigma=\d\stars\sigma=0$ 
is flat to first order at~$0\in\bbR^7$.

Now, if~$\psi$ is any $3$-form on~$\bbR^7$ 
that vanishes at the origin~$0\in\bbR^7$, 
then, because~$\L^3_+(V^*)$ is an open set in~$\L^3(V^*)$,
the $3$-form~$\sigma = \phi + \psi$ is a definite $3$-form 
on some open neighborhood of~$0\in\bbR^7$.  Since~$\d\sigma
= \d\psi$, and since, for any~$4$-form~$\Psi\in\L^4(V^*)$,
there exists a $3$-form~$\psi$ on~$\bbR^7$ that vanishes
at~$0\in\bbR^7$ and that satisfies~$(\d\psi)_0 = \Psi$,
it follows that the condition~$\d\sigma = 0$,
i.e.,~$\tau_0 = \tau_1 = \tau_3 = 0$, imposes $35$
independent linear conditions on the intrinsic torsion
of~$\sigma$.  Since these conditions must define some 
\Gtwo-invariant subspace of the torsion representation
$\Vs_{0,0}{\oplus}\Vs_{1,0}{\oplus}\Vs_{0,1}{\oplus}\Vs_{2,0}$,
it follows by dimension count that it is the subspace
$\Vs_{0,0}{\oplus}\Vs_{1,0}{\oplus}\Vs_{2,0}$.

Similarly, since~$\L^4_+(W^*)$ is an open subset of~$\L^4(W^*)$
and since~$\Ss:\L^3_+(W^*)\to\L^4_+(W^*)$ is a smooth double
covering, it follows that if~$\psi$ is any smooth $4$-form 
vanishing at the origin~$0\in\bbR^7$, then there is an open
neighborhood~$U$ of~$0\in\bbR^7$ on which there exists a
definite form~$\sigma$ such that~$\sigma_0 = \phi$ and
$\stars\sigma = {\ast_\phi}\phi + \psi$.  Moreover, if~$\Psi$
is any $5$-form in $\L^5(V^*)$, then there exists a smooth
$4$-form~$\psi$ vanishing at~$0\in\bbR^7$ such that~$(\d\psi)_0
=\Psi$.  The corresponding definite $3$-form~$\sigma$ 
will then satisfy~$\d\stars\sigma = \d\psi$, 
so that~$(\d\stars\sigma)_0 = \Psi$.  It follows that the
condition~$\d\stars\sigma=0$, i.e., $\tau_1 = \tau_2 = 0$,
must be $21$ independent linear equations on the intrinsic
torsion of~$\sigma$.  Since these conditions must define some 
\Gtwo-invariant subspace of the torsion representation
$\Vs_{0,0}{\oplus}\Vs_{1,0}{\oplus}\Vs_{0,1}{\oplus}\Vs_{2,0}$,
it follows by dimension count that it is the subspace
$\Vs_{1,0}{\oplus}\Vs_{0,1}$.  

Thus, the conditions~$\d\sigma=0$ and~$\d\stars\sigma=0$
together imply that all of the intrinsic torsion of~$\sigma$
vanishes, i.e., that~$\sigma$ is flat to first order at each point.
\end{proof}

\begin{remark}[Fern\'andez and Gray's theorem on vector cross products]
\label{rem: grayresult}
Proposition~\ref{prop: 1flatg2} implies the 1982 result 
of Fern\'andez and Gray~\cite{MR84e:53056}
that a vector cross product~$\times:TM{\times}TM\to TM$ 
that is compatible with a Riemannian metric~$g$ on~$M$ is $g$-parallel
if and only if the corresponding~$3$-form is closed and coclosed
(with respect to~$g$).

The essential difference between Proposition~\ref{prop: 1flatg2}
and their result is that they assume a specific metric~$g$ 
and vector cross product to be given, whereas Proposition~\ref{prop: 1flatg2}
starts with a definite $3$-form~$\sigma$ and constructs a specific metric
associated to~$\sigma$.
\end{remark}

\section{Frame Bundle Calculations} 

\subsection{The associated Levi-Civita connection}
Let~$\sigma\in\Omega^3_+(M)$ be a \Gtwo-structure with
associated \Gtwo-bundle~$F_\sigma\subset\cF$.  This
bundle can be canonically enlarged to an oriented
orthonormal frame bundle~$\Fs_\sigma = F_\sigma\cdot SO(7)\subset\cF$
and this larger bundle will be referred to as the 
associated metric frame bundle of~$\sigma$. 

Now $\pi\colon\Fs_\sigma\to M$ has a tautological $V$-valued 
$1$-form $\eomega$ defined by requiring that~$\eomega(v) 
= u(\pi_*(v))$ for all~$v\in T_u\Fs$.  
It may help the reader to think of $\eomega$ as expanded in the 
basis~$e_i$ in the form~$\eomega = \eomega_1\,e_1+\cdots+\eomega_7\,e_7$ 
and then think of~$\eomega$ as a column of height~$7$, i.e., 
$\eomega=(\eomega_i)$.  

The Levi-Civita connection is then represented on $\Fs_\sigma$ 
as a $1$-form~$\epsi$ on~$\Fs_\sigma$ taking
values in $\euso(7)$, i.e., the $7$-by-$7$ skew-symmetric matrices.
As such,~$\epsi = (\epsi_{ij})$ where~$\epsi_{ij}=-\epsi_{ji}$.

The defining property of~$\epsi$ is that it satisfies
the \emph{first structure equation} of Cartan:
\begin{equation}\label{eq: 1st str eq}
   \d\,\eomega = -\epsi\w\eomega.
\end{equation}
In indices (i.e., components) this matrix equation becomes 
the system of equations~$\d\eomega_i=-\epsi_{ij}\w\eomega_j$.

The curvature of this connection is represented by the 
$2$-form~$\ePsi = \d\epsi + \epsi \w \epsi$.  It satisfies
the \emph{first Bianchi identity}
\begin{equation}\label{eq: 1st Bianchi metric}
\ePsi \w \eomega = 0
\end{equation}
and has the indicial
expression
\begin{equation}
\ePsi_{ij} = \d\epsi_{ij} + \epsi_{ik}\w\epsi_{kj}
= {\ts\frac12}\eR_{ijkl}\,\eomega_k\w\eomega_l\,.
\end{equation}

\subsubsection{The natural connection and intrinsic torsion on~$F_\sigma$}
To save writing, I will denote the pullbacks
of~$\eomega$ and~$\epsi$ to~$F_\sigma$ by the same letters,
trusting the reader to keep in mind where various equations
are taking place.  

The pullback of $\epsi$ to $F_\sigma$ will not generally have 
values in $\eug_2\subset\euso(7)$.  However, keeping
in mind the canonical decomposition~$\euso(7) = \eug_2\oplus [V]$,
there is a unique decomposition of the form 
\begin{equation}\label{eq: 1st str eq -- tor}
 \epsi = \etheta + 2[\etau]
\end{equation}
where $\etheta$ takes values in~$\eug_2$ and~$\etau$ takes values 
in~$V$.  (The coefficient~$2$ simplifies subsequent formulas.)

Then~$\etheta$ is a connection $1$-form on~$F_\sigma$ and defines
what will be referred to as the \emph{natural} connection associated
to the \Gtwo-structure~$\sigma$.  This connection will not be
torsion-free (and hence is not the Levi-Civita connection) 
unless~$\etau$ vanishes identically. 

\subsection{General $G$-structure torsion}\label{ssec: bundle torsion}
This construction of a natural connection for a \Gtwo-structure~$\sigma$
is an instance of a general construction valid for any~$G\subset\Or(n)$.

Letting~$\eug\subset\euso(n)$ denote the Lie algebra of~$G$, 
there is a unique $G$-equivariant splitting
$\euso(n) = \eug\oplus\eug^\perp$ obtained by using 
the standard $\Or(n)$-invariant inner product on~$\euso(n)$.  

For any $G$-structure~$\pi:F\to M$, one has the associated
orthonormal frame bundle~$\Fs = F\cdot\Or(n)$.  One can then pull
back the Levi-Civita connection~$\epsi$ on~$\Fs$ to~$F$ and
decompose it uniquely in the form~$\epsi = \etheta + \tau$
where~$\etheta$ takes values in~$\eug$ and~$\tau$ takes values in
$\eug^\perp\simeq \euso(n)/\eug$.  The $1$-form~$\etheta$ defines
a natural connection on~$F$ (one that is the pullback to~$F$ of
a metric-compatible connection, generally with torsion, on~$\Fs$).
The $1$-form~$\etau$ represents a section~$T$ of the associated 
\emph{torsion bundle}~$F\times_\rho (\eug^\perp\otimes\bbR^n)$,
where~$\rho:G\to\End\bigl(\eug^\perp\otimes\bbR^n\bigr)$ is the
tensor product of the two obvious representations.

It is a general result (essentially due to \'E. Cartan) 
that all of the pointwise first-order diffeomorphism 
invariants of a $G$-structure~$F\subset\cF$ that are polynomial 
in the derivatives of the corresponding defining section~$\sigma$
of the bundle~$\cF/G$ are expressible as polynomials in the section~$T$.

Moreover, for~$k\ge 2$, all of the pointwise $k$-th order diffeomorphism 
invariants of a $G$-structure~$F\subset\cF$ that are polynomial 
in the first $k$ derivatives of the corresponding defining 
section~$\sigma$ of the bundle~$\cF/G$ are expressible as polynomials 
in the section~$T$, its first~$k{-}1$ covariant derivatives with
respect to the connection~$\etheta$, the curvature of~$\etheta$,
and its first $k{-}2$ covariant derivatives (with respect to~$\etheta$).

Consequently, for each~$k\ge1$, the polynomial pointwise invariants 
of order~$k$ are polynomials in a canonically defined section of 
a vector bundle of the form
\begin{equation*}
F\times_{\rho_1\times\cdots\times\rho_k}
  \bigl(V_1(\eug)\oplus\cdots\oplus V_k(\eug)\bigr)
\end{equation*} 
where~$V_k(\eug)$ is 
the unique $G$-representation that satisfies
\begin{equation}\label{eq: kth order invariants}
\bigl(\eugl(n,\bbR)/\eug\bigr) \otimes \Ss^k(\bbR^n)
= V_k(\eug) \oplus \bigl(\bbR^n\otimes \Ss^{k+1}(\bbR^n)\bigr).
\end{equation}

In the familiar case in which~$\eug=\euso(n)$, 
the first torsion space~$V_1\bigl(\euso(n)\bigr)$ vanishes
(this is simply the fundamental lemma of Riemannian geometry)
and one has the result (due to Cartan and Weyl) that all 
of the pointwise invariants of a metric can be expressed in 
terms of the Riemann curvature tensor and its covariant 
derivatives with respect to the Levi-Civita connection.

\begin{remark}[Canonical connections]\label{rem: natural vs canonical}
The use of the term `natural' with regard to the connection~$\etheta$
on the $G$-structure~$F$ should not be construed to mean that this
is the only `canonical' connection on~$M$ that is compatible 
with~$F$.  In many cases, this is only one of a family of possible
`canonical' connections that can be defined in terms of the first-order
invariants of the~$G$-structure~$F$ and that are preserved under
equivalence of $G$-structures.  

For example, if the $G$-modules~$V_1(\eug)$ 
and~$\eug\otimes\bbR^n$ have common constituents, so that the
space~$\Hom^G\bigl(V_1(\eug),\eug\otimes\bbR^n\bigr)$ of $G$-equivariant
homomorphisms between the two spaces has dimension~$r>0$, there
will be an $r$-parameter family of ways of modifying~$\etheta$,
by adding a~$\eug$-valued~$1$-form whose coefficients are linear
in the torsion functions, in such a way that the resulting modification
defines a connection on~$M$ compatible with the $G$-structure~$F$.
Each element in this $r$-parameter family of connections can be 
regarded as canonical in the sense that equivalence of~$G$-structures
will induce isomorphisms between the corresponding connections in
the $r$-parameter family.

Of course, there is no \emph{a priori} reason to 
consider only connection modifications that are linear 
in the torsion functions; for example, any $G$-equivariant polynomial 
mapping~$V_1(\eug)\to\eug\otimes\bbR^n$ could be used to 
define such a modification of~$\etheta$.  However, these
`higher' modifications do not often arise in practice.

Depending on the intended use, it could
well be that one of these other connections (rather than the one
being called `natural' in the present article) is better suited for
expressing identities of one kind or another.
\end{remark}

\subsection{\Gtwo-specific calculations}
In the specific case of~$\Gtwo\subset\SO(7)$, one finds, 
as has already been remarked,
\begin{equation}\label{eq: V1 of G2}
V_1(\eug_2) \simeq \Vs_{0,0}\oplus\Vs_{1,0}\oplus\Vs_{0,1}\oplus\Vs_{2,0},
\end{equation}
while~$V_2(\eug_2)$, which has dimension~$392$, has the decomposition
\begin{equation}\label{eq: V2 of G2}
V_2(\eug_2) \simeq  \Vs_{0,0}\oplus 2\Vs_{1,0}\oplus
                    \Vs_{0,1}\oplus 3\Vs_{2,0}\oplus
                    2\Vs_{1,1}\oplus \Vs_{0,2}\oplus \Vs_{3,0}\,.
\end{equation}
Naturally, this latter space has~$V_2\bigl(\euso(7)\bigr)$,
i.e., the curvature tensors of metrics in dimension~$7$, as a quotient.
For comparison, note that, as \Gtwo-modules:
\begin{equation}\label{eq: V2 of SO7}
V_2\bigl(\euso(7)\bigr)
            \simeq  \Vs_{0,0}\oplus 2\Vs_{2,0}\oplus
                    \Vs_{1,1}\oplus \Vs_{0,2}\,.
\end{equation}
The Ricci tensor takes values in a subspace isomorphic 
to~$\Vs_{0,0}\oplus\Vs_{2,0}$ while the remainder
represents the Weyl tensor.

\begin{remark}[Canonical $\Gtwo$-connections]\label{rem: canonical g2 conns}
Since~$\eug_2\otimes\Vs_{1,0}=\Vs_{1,0}\oplus\Vs_{2,0}\oplus\Vs_{1,1}$
shares two $\Gtwo$-irreducible modules with $V_1(\eug)$, it follows from
Remark~\ref{rem: natural vs canonical} that there is actually a $2$-parameter
family of canonical connections associated to any~$\Gtwo$-structure~$\sigma$. 
Each element in this family is compatible with~$\sigma$ (in the sense
that~$\sigma$ is parallel under the corresponding parallel translation).  
Since the common constituents~$\Vs_{1,0}$ and~$\Vs_{2,0}$ correspond to the 
torsion forms~$\tau_1$ and~$\tau_3$, respectively, it follows that the
entire two-parameter family of canonical connections collapses to a single
connection if and only if the $\Gtwo$-structure~$\sigma$ 
satisfies~$\tau_1=\tau_3=0$.
In this case, differentiating the equations~\eqref{eq: dsigma delsigma} 
shows that $\tau_0\tau_2=0$ and~$\d\tau_0=0$.  In particular, when~$M$ is connected, it follows that~$\tau_0$ is constant.  If~$\tau_0=0$, then
the $\Gtwo$-structure is closed.  If~$\tau_0\not=0$, then~$\tau_2=0$
and one has the equation~$\d\sigma = \tau_0\,\stars\sigma$, which 
is the defining equation for the so-called `nearly $\Gtwo$-manifolds'.  

Thus, the family of canonical $\Gtwo$-connections associated to a
$\Gtwo$-structure~$\sigma$ collapses to a single $\Gtwo$-connection 
if and only if either~$\sigma$ is closed or it defines a nearly 
$\Gtwo$-manifold.
\end{remark}

\subsection{The second structure equations}
It is helpful to make the following observation:  
The identities~\eqref{eq: box and angle identities}  
imply that the 2-form $2[\etau]\w[\etau] + [\,[\etau]\w\etau\,]$ 
takes values in $\eug_2$.  This motivates the definitions
\begin{align}\label{eq: Dtau and Dtheta}
D\etau &= \d\etau + \etheta\w\etau - [\etau]\w\etau\\
D\etheta &= \d\etheta + \etheta\w\etheta + 4 [\etau]\w[\etau] + 
2[\,[\etau]\w\etau\,],
\end{align}
for, with these definitions, 
$D\etheta$ takes values in~$\eug_2$. Moreover
\begin{equation}
\ePsi = d(\etheta + 2[\etau]) 
+ (\etheta + 2[\etau])\w(\etheta + 2[\etau]) = 
D\etheta + 2[D\etau]
\end{equation}
so that the first Bianchi identity takes the form
\begin{equation}
(D\etheta + 2[D\etau])\w\eomega = 0.
\end{equation}

\begin{remark}[Covariant differentials]\label{rem: cov diffs}
The decisive advantage of using the forms~$D\etau$ and~$D\etheta$
to express the curvature tensor is that these forms do \emph{not}
contain all of the information about the second order invariants
of the underlying~\Gtwo-structure~$\sigma$ although they do
contain enough information to recover the Riemann curvature
tensor of the underlying metric.
\end{remark}

\subsection{Indicial calculations}
The indicial expression of~\eqref{eq: 1st str eq}
in terms of~\eqref{eq: 1st str eq -- tor} is
\begin{equation}\label{eq: indicial 1st str eq}
\begin{aligned}
\d\eomega_i = -\etheta_{ij}\w\eomega_{j} - 2\eps{ijk}\etau_k\w\eomega_j\,.
\end{aligned}
\end{equation}

Denote $\pi^*(\sigma)$ by $\esigma$ and, with a slight abuse 
of notation, denote~$\pi^*(\stars\sigma)$ by $\star\esigma$. 
Then
\begin{align}
\esigma &= 
{\ts\frac{1}{6}}\eps{ijk}\,\eomega_i\w\eomega_j\w\eomega_k\\
\star\esigma &= 
{\ts\frac{1}{24}}\eps{ijkl}\,
\eomega_i\w\eomega_j\w\eomega_k\w\eomega_l\,.
\end{align}
These give rise, via~\eqref{eq: indicial 1st str eq} and the 
$\eps{}$-identities, to the formulae
\begin{equation}
\begin{aligned}
\d\,\esigma &= 
\eps{ijkl}\,\etau_i\w\eomega_j\w\eomega_k\w\eomega_l\\
\d\,{\star}\esigma &= 
-\left(\etau_p\w\eomega_p\right)\>\w\>
\left(\eps{ijk}\,\eomega_i\w\eomega_j\w\eomega_k\right)\\
&= -6\left(\etau_p\w\eomega_p\right) \w \esigma
\end{aligned}
\end{equation}

\subsubsection{Torsion decomposition.}
There are unique functions~$\eT_{ij}$ on~$F_\sigma$ so that
\begin{equation}
\etau_i = \eT_{ij}\,\eomega_j\,.
\end{equation}

These functions can be used to express the 
intrinsic torsion forms in indicial form:
\begin{equation}\label{eq: indicial torsion forms}
\begin{aligned}
\pi^*(\tau_0) &= {\ts\frac{24}7}\,\eT_{ii}\,,\\
\pi^*(\tau_1) &= \eps{ijk}\eT_{ij}\,\,\omega_k\,,\\
\pi^*(\tau_2) &= 4\,\eT_{ij}\,\omega_i\w\omega_j 
                  - \eps{ijkl}\eT_{ij}\,\omega_k\w\omega_l\,,\\
\pi^*(\tau_3) &=  -{\ts\frac32}\eps{ikl}\,(\eT_{ij}+\eT_{ji})
                  \,\eomega_j\w\eomega_k\w\eomega_l
                   +{\ts\frac{18}7}\,\eT_{ii}\,\esigma\,.
\end{aligned}
\end{equation}
(In these formulae, one sums over repeated indices in
any term.)

\subsubsection{Curvature identities}
The covariant differentials can be expressed in indices as
\begin{align}\label{eq: Tijk def}
D\etau &= (D\etau_i) = 
({\ts\frac12}\eT_{ijk}\,\eomega_j\w\eomega_k)\\
\label{eq: Sijkl def}
D\etheta &= (D\etheta_{ij}) = 
({\ts\frac12}\eS_{ijkl}\,\eomega_j\w\eomega_k)
\end{align}
where each of $\eT$ and $\eS$ are skew-symmetric in their last two 
indices, $\eS$ is skew-symmetric in its first two indices and, since 
$D\etheta$ takes values in $\eug_2$, the functions $\eS$ also satisfy
$$\eps{ijm}\eS_{ijkl}=0$$ for all $m$, $k$, and $l$.  

Since~$\ePsi = D\etheta + 2\bigl[D\etau\bigr]$,
the Riemann curvature functions are expressed as
\begin{equation}\label{eq: R is S + 2 T}
\eR_{ijkl} = \eS_{ijkl} + 2\,\eps{ijp}\eT_{pkl}\,,
\end{equation}
so that the first Bianchi identity becomes
\begin{equation}\label{eq: indicial 1st Bianchi}
\eS_{ijkl}+\eS_{iljk}+\eS_{iklj}+2\,\eps{ijp}\eT_{pkl}+
2\,\eps{ilp}\eT_{pjk}+2\,\eps{ikp}\eT_{plj}=0.
\end{equation}

The identities~\eqref{eq: indicial 1st Bianchi} 
impose $28$ linear conditions on the $\eT_{ijk}$ 
alone.  Perhaps the easiest way to derive these $28$ 
conditions is to expand the identities 
\begin{equation}\label{eq: dds = ddss = 0}
\d(\d(\esigma))=\d(\d(\star\esigma))=0
\end{equation} 
and use the structure equations~\eqref{eq: indicial 1st str eq}
together with the 
definitions~\eqref{eq: Dtau and Dtheta},~\eqref{eq: Tijk def},
and~\eqref{eq: Sijkl def}. This will be left as an 
exercise for the reader.  The result is that the 
conditions~\eqref{eq: dds = ddss = 0} are equivalent 
to the following equations (some of which are redundant):
\begin{equation}\label{eq: indicial Tijk relns}
\begin{aligned}
0 &= \eT_{iij}\,,\\
0 &= \eps{ipq}\eT_{jpq} - \eps{jpq}\eT_{ipq}\,,\\ 
0 &= \eps{ipq}\eT_{pqj} - \eps{jpq}\eT_{pqi}\,.
\end{aligned}
\end{equation}
This implies that the function~$(\eT_{ijk})$,
which nominally takes values in a \Gtwo-module of the 
form
\begin{equation}
\begin{aligned}
\Vs_{1,0}\otimes\L^2\bigl(\Vs_{1,0}\bigr)
&=\Vs_{1,0}\otimes\bigl(\Vs_{1,0}\oplus\Vs_{0,1}\bigr)\\
&= \Vs_{0,0}\oplus 2\Vs_{1,0}\oplus\Vs_{0,1}
      \oplus 2\Vs_{2,0}\oplus\Vs_{1,1}\,,
\end{aligned}
\end{equation}
actually takes values in a submodule of the form
\begin{equation}
       \Vs_{0,0}\oplus 2\Vs_{2,0}\oplus\Vs_{1,1}\,.
\end{equation}

\subsubsection{The Ricci identity}
It was Bonan~\cite{MR33 4855} who first observed that the Bianchi
identities imply that a \Gtwo-structure with vanishing torsion 
must necessarily have vanishing Ricci tensor.
On general abstract grounds, it then follows that the Bianchi 
identities~\eqref{eq: indicial 1st Bianchi} must allow one
to express the Ricci curvature in terms of the~$\eT_{ijk}$.
Indeed, by combining the first Bianchi identities via the 
$\eps{}$-identities (another exercise for the reader), one 
derives the following expression for 
the Ricci curvature components $\eR_{ij} = \eR_{kikj}$:
\begin{equation}\label{eq: indicial Ricci as Tijk}
\eR_{ij} = 6\,\eps{pqi}\eT_{pqj}\,.
\end{equation}

This allows one to express the Ricci curvature
directly in terms of the four torsion forms and their
exterior derivatives. The resulting formula for the scalar 
curvature of the underlying metric~$g_\sigma$ is
\begin{equation}\label{eq: scal curv form}
\Scal(g_\sigma) = 12\,\delta\tau_1 + {\ts\frac{21}{8}}\,{\tau_0}^2 
     + 30\,|\tau_1|^2-{\ts\frac12}\,|\tau_2|^2 -{\ts\frac12}\,|\tau_3|^2.
\end{equation}

The full Ricci tensor is somewhat more 
complicated, but can be expressed as follows:  

First, define a \Gtwo-invariant quadratic 
pairing~$\Qs:\L^3(T^*)\times\L^3(T^*)\to\L^3(T^*)$
by the following recipe:  Choose a local basis~$e_1,\ldots,e_7$ of 
orthonormal vector fields such that~$\sigma(e_i,e_j,e_k)=\eps{ijk}$
(such a basis is often called a \Gtwo-frame field).  Then, for
$\alpha,\beta\in\Omega^3(M)$ set
\begin{equation}\label{eq: Q def}
\Qs(\alpha,\beta) 
= \stars\bigl[\eps{ijkl}\,\bigl((e_i\w e_j)\lhk\stars\alpha\bigr)
                \w \bigl((e_k\w e_l)\lhk\stars\beta\bigr)\bigr]\,.
\end{equation}
The resulting mapping~$\Qs$ does not depend on the choice
of local~\Gtwo-frame field.  With this definition
(and keeping in mind the definition~\eqref{eq: js def} of~$\js$) 
one finds
\begin{equation}\label{eq: Ric formula}
\begin{aligned}
\Ric(g_\sigma) &= 
  -\biggl({\ts\frac32}\,\delta\tau_1 
  -{\ts\frac38}\,{\tau_0}^2 + 15\,|\tau_1|^2
  -{\ts\frac14}\,|\tau_2|^2 +{\ts\frac12}\,|\tau_3|^2\biggl)\,g_\sigma\\
  \ & \qquad
    + \js\Bigl(\ - {\ts\frac54}\,\d\bigl(\stars(\tau_1\w\stars\sigma)\bigr)
               - {\ts\frac14}\,\d\tau_2  + {\ts\frac14}\,\stars\d\tau_3\\
     & \qquad\qquad\ \ 
               + {\ts\frac52}\,\tau_1 \w \stars(\tau_1\w\stars\sigma)
               - {\ts\frac18}\,\tau_0\tau_3
               + {\ts\frac14}\,\tau_1\w\tau_2 \\
     & \qquad\qquad\quad 
               + {\ts\frac34}\,\stars(\tau_1\w\tau_3)
               + {\ts\frac18}\,\stars(\tau_2\w\tau_2)
               + {\ts\frac1{64}}\,\Qs(\tau_3,\tau_3)\ \Bigr)\,.
\end{aligned}
\end{equation}
While a formula in this generality is not of much practical
use, when one goes to investigate special classes of 
\Gtwo-structures, this formula can simplify considerably, as
will be seen.

Formulae essentially equivalent to a special case of the 
formulae~\eqref{eq: scal curv form} and~\eqref{eq: Ric formula}
were found in~\cite{math.DG/0102142,math.DG/0112201}, where those
authors considered what they called `integrable $\Gtwo$-structures',
which, in the terminology of this article, means 
$\Gtwo$-structures~$\sigma$ satisfying~$\tau_2=0$. 

\begin{remark}[General identities]\label{rem: general identities}
It is perhaps worth remarking on why the 
identities~\eqref{eq: scal curv form} and~\eqref{eq: Ric formula}
could be expected to have the form that they do.  

In the first place, one knows that the scalar curvature
must be expressible in a \Gtwo-invariant manner as a
sum of a linear expression in the second order invariants, i.e.,
a section of a vector bundle modeled on~$V_2(\eug_2)$, and
an expression in the first order invariants, i.e.,
the torsion forms, that is at most quadratic.  
A glance at~\eqref{eq: V2 of G2} shows that there is only 
one trivial summand in the representation~$V_2(\eug_2)$ 
and hence there is essentially only one possible second
order term up to a universal constant multiple.  Since~$\delta\tau_1$
is a scalar second order invariant, it must represent this
copy of~$\Vs_{0,0}$ in~$V_2(\eug_2)$.  As for the first
order terms, since~$V_1(\eug_2)$ consists of four mutually
inequivalent~\Gtwo-modules, the space of \Gtwo-invariant quadratic
forms  on this space has dimension~$4$ and must be represented
by the square norms of the four torsion forms.  Thus, a
formula of the form~\eqref{eq: scal curv form} was inevitable; 
it was just a matter of determining the numerical coefficients,
which was done with the aid of~\textsc{Maple}.

The argument for the form of~\eqref{eq: Ric formula} is 
quite similar.  Since the scalar curvature has already been
determined, it is a question of writing down a formula for
the trace-free part of the Ricci tensor, i.e., finding 
linear terms in~$V_2(\eug_2)$ and quadratic terms 
in~$V_1(\eug_2)$ that take values in the \Gtwo-module~$\Vs_{2,0}$.
Again, a glance at~\eqref{eq: V2 of G2} shows that
there are at most three possible second order terms and it
is not difficult to see that the three second order terms
that take values in $\Vs_{2,0}$ found by taking
derivatives of~$\tau_1$, $\tau_2$, and $\tau_3$ and projecting
into a suitable~$\Vs_{2,0}$ representation are, in fact, 
independent and generate the three copies of~$\Vs_{2,0}$
that appear in~$V_2(\eug_2)$.  On the other hand, using 
representation theory to compute the second symmetric power
of~$V_1(\eug_2)$ shows that there exist eight copies of~$\Vs_{2,0}$
in this symmetric power.  Of those eight copies, five are
computable via wedge product and appear in the formula for
Ricci.  Of the remaining three, one bilinear in $\tau_2$
and~$\tau_3$ and the other two quadratic in~$\tau_3$, only
one of the terms quadratic in~$\tau_3$ actually makes an
appearance.  The rest is just a matter of determining constants.
\end{remark}

\subsection{Closed \Gtwo-structures}\label{ssec: closed G2 structures}
Now, consider the case of a 
\emph{closed}~$\sigma\in\Omega^3_+(M)$, i.e.,~$\d \sigma = 0$.
In this case, by Proposition~\ref{prop: ext der formulae}, 
it follows that
\begin{equation}\label{eq: def tau2}
\d\stars\sigma = \tau_2\w\sigma,
\end{equation}
where~$\tau_2$ lies in~$\Omega^2_{14}(M,\sigma)$.  
In particular,
\begin{equation}\label{eq: tau2 no 7}
\tau_2\w\stars\sigma=0.
\end{equation}
Taking the exterior derivative of~\eqref{eq: def tau2} yields
\begin{equation}\label{eq: dtau2 no 7}
0 = \d\tau_2 \w \sigma\,,
\end{equation}
implying that~$\d\tau_2$ has no component in~$\Omega^3_7(M,\sigma)$.  
Differentiating~\eqref{eq: tau2 no 7} yields
\begin{equation}\label{eq: identify 1 comp dtau2}
\begin{aligned}
0 
 = \d\bigl(\tau_2\w\stars\sigma\bigr) 
&= \d\tau_2 \w\stars\sigma + \tau_2\w \d\stars\sigma\\
&= \d\tau_2 \w\stars\sigma + \tau_2\w \tau_2\w\sigma
 = \d\tau_2 \w\stars\sigma - |\tau_2|^2\,\stars1\,.
\end{aligned}
\end{equation}
Thus, from~\eqref{eq: dtau2 no 7} and~\eqref{eq: identify 1 comp dtau2}
it follows that there exists a $\gamma\in\Omega^3_{27}(M,\sigma)$ so that
\begin{equation}\label{eq: 1 part dtau2}
\d\tau_2 = {\ts\frac17}|\tau_2|^2\,\sigma + \gamma.
\end{equation}

In summary, formulae~\eqref{eq: scal curv form} 
and~\eqref{eq: Ric formula} can be simplified in this case to
\begin{equation}\label{eq: Scal when closed}
\Scal(g_\sigma) = -{\ts\frac12}\,|\tau_2|^2.
\end{equation}
and
\begin{equation}\label{eq: Ric when closed}
\Ric(g_\sigma) = {\ts\frac14}\,|\tau_2|^2 \,g_\sigma
    -{\ts\frac14}\,\js\left(\,\d\tau_2 
               -{\ts\frac12}\,\stars(\tau_2\w\tau_2)\,\right).
\end{equation}

\begin{remark}[Differential invariants of closed \Gtwo-structures]
\label{rem: diff inv closed structures}
Just as one can compute the dimension of the space of $k$-jets
of $G$-structures as in \S\ref{ssec: bundle torsion}, one can
compute the dimension of the space of $k$-jets of $G$-structures
satisfying some set of differential equations.  In the case
of closed \Gtwo-structures, denote the module of $k$-th order
differential invariants by~$V'_k(\eug_2)\subset V_k(\eug_2)$.
One finds, for example, that
\begin{equation}\label{eq: V1p and V2p}
V'_1(\eug_2)\simeq \Vs_{0,1}\qquad\qquad 
V'_2(\eug_2)\simeq \Vs_{2,0}\oplus \Vs_{1,1} \oplus \Vs_{0,2}.
\end{equation}
This implies, on abstract grounds, that
the scalar curvature of the underlying metric of a closed 
\Gtwo-structure must be expressed in terms of the first order
invariants (since there is no $\Vs_{0,0}$ component in~$V'_2(\eug_2)$)
and that the full Ricci tensor can be expressed in terms 
of~$\tau_2$ and~$\d\tau_2$. Thus, the form of~\eqref{eq: Scal when closed} 
and~\eqref{eq: Ric when closed} could have been anticipated,
if not the numerical coefficients.

Of course, it is easy to `write down' the general closed 
\Gtwo-structure locally:  If~$\beta\in\Omega^2(\bbR^7)$ is
a (smooth) $2$-form that vanishes to second order at~$0\in\bbR^7$,
then the $3$-form~$\sigma = \phi + \d\beta$ will equal~$\phi$
at~$0$ and hence will be a closed, definite $3$-form on some
open neighborhood of~$0\in\bbR^7$.  Conversely, if~$\sigma$ 
is a closed \Gtwo-structure on a manifold~$M^7$, then any point~$p\in M$ 
has an open neighborhood~$U$ on which there exists a $p$-centered 
coordinate chart~$x:U\to\bbR^7$ such that~$\sigma_{U}= x^*(\phi+\d\beta)$ where~$\beta\in\Omega^2(\bbR^7)$ is a $2$-form that vanishes 
to second order at~$0\in\bbR^7$.

In a sense that it is possible to make precise using Cartan's
notion of the generality of the space of solutions of a system
of~\textsc{pde}, one can develop this discussion further to
show that the general closed \Gtwo-structure modulo diffeomorphism
depends on $8$ functions of seven variables.
\end{remark}

An immediate consequence of~\eqref{eq: Scal when closed} 
is the following:

\begin{corollary}\label{cor: closed => sc <= 0}
For any closed \Gtwo-structure~$\sigma\in\Omega^3_+(M)$, 
the scalar curvature of the underlying metric is non-positive 
and vanishes identically 
if and only if the entire Ricci tensor of the underlying metric
vanishes.  Equivalently, the scalar curvature vanishes identically 
if and only if $\sigma$ satisfies~$\d\sigma=\d\stars\sigma=0$.
\end{corollary}

Using the formulae for~$\is$ and~$\js$, the
formula~\eqref{eq: Ric when closed}
can be rewritten as
\begin{equation}\label{eq: dtau as Ric when closed}
\d\tau_2 = {\ts\frac3{14}}\,|\tau_2|^2\,\sigma
+{\ts\frac12}\stars(\tau_2\w\tau_2) 
-{\ts\frac12}\is\bigl(\Ric^0(g_\sigma)\bigr),
\end{equation}
where~$\Ric^0(g_\sigma)$ is the traceless Ricci tensor of~$g_\sigma$.

\begin{corollary}\label{cor: closed Einstein}
A closed \Gtwo-structure~$\sigma\in\Omega^3_+(M)$ has
an Einstein underlying metric
if and only if it satisfies~$\d\stars\sigma = \tau_2\w\sigma$
where $\d\tau_2 = {\ts\frac3{14}}\,|\tau_2|^2\,\sigma
+{\ts\frac12}\stars(\tau_2\w\tau_2)$.
\end{corollary}

\begin{remark}[Nonexistence of compact Einstein examples]
\label{rem: closed Einstein?}
I do not know whether there exist any closed
\Gtwo-structures that are Einstein but not Ricci-flat, even
local (i.e., incomplete) ones. 

After Version~$1.0$ of the present article was posted to 
the arXiv, Cleyton and Ivanov~\cite{math.DG/0306362} gave
an argument (based on a comparison of the Ricci curvatures
of the Levi-Civita connection and the canonical connection
of the underlying $\Gtwo$-structure) showing that no
\emph{compact} $7$-manifold can support a closed 
$\Gtwo$-structure~$\sigma$ whose underlying metric~$g_\sigma$ 
is Einstein unless $\sigma$ is also coclosed, i.e., $\d\stars\sigma=0$.
Their argument is rather involved, but 
Corollary~\ref{cor: closed Einstein} yields 
a simple proof:

Suppose that~$\sigma\in\Omega^3_+(M)$ is a closed $\Gtwo$-structure
whose underlying metric~$g_\sigma$ is Einstein.  Then, by 
Corollary~\ref{cor: closed Einstein}, it follows that 
$\d\stars\sigma = \tau_2\w\sigma$
where $\d\tau_2 = {\ts\frac3{14}}\,|\tau_2|^2\,\sigma
+{\ts\frac12}\stars(\tau_2\w\tau_2)$.  Now, using this formula
together with the formula~\eqref{eq: G2 quartic}, one finds
\begin{equation}\label{eq: d tau cubed}
\begin{aligned}
\d\bigl({\ts\frac13}{\tau_2}^3\bigr)
&= {\tau_2}^2 \w \d \tau_2 
 = {\tau_2}^2 \w\bigl({\ts\frac3{14}}\,|\tau_2|^2\,\sigma
+{\ts\frac12}\stars(\tau_2\w\tau_2)\bigr)\\
&= -{\ts\frac3{14}}\,|\tau_2|^4\stars 1 
   +{\ts\frac12} |\tau_2\w\tau_2|^2\,\stars 1
= {\ts\frac27} |\tau_2|^4\stars 1.
\end{aligned}
\end{equation}
Now, suppose that~$M$ were compact.  Integrating both ends of
\eqref{eq: d tau cubed} over~$M$ and applying Stokes' theorem yields
\begin{equation}
0 = \int_M \d\bigl({\ts\frac13}{\tau_2}^3\bigr)
= \int_M {\ts\frac27} |\tau_2|^4\stars 1\,,
\end{equation}
implying that $\tau_2$ must vanish identically, as was to be shown.
\end{remark}

In view of \eqref{eq: dtau as Ric when closed}, 
this nonexistence can be seen as a special case of a
general result about pinching of Ricci curvature:

\begin{corollary}\label{cor: pinched Ricci}
Suppose that~$\sigma\in\Omega^3_+(M)$ is a closed $\Gtwo$-structure
on a compact $7$-manifold~$M$ that satisfies the pinching condition
\begin{equation}\label{eq: Ricci pinch ineq}
\bigl|\Ric^0(g_\sigma)\bigr|^2 \le {\ts\frac{4}{21}}\,C\,\Scal(g_\sigma)^2.
\end{equation}
for some constant~$C\le 1$.  If~$C<1$, then $\sigma$ is also 
coclosed.  If~$C=1$, then equality must hold 
in~\eqref{eq: Ricci pinch ineq} everywhere on~$M$.  Moreover, 
in this case, the identity
\begin{equation}
\is\bigl(\Ric^0(g_\sigma)\bigr) 
= {\ts\frac23}\bigl(\stars(\tau_2\w\tau_2)+{\ts\frac17}|\tau_2|^2\,\sigma\bigr)
\end{equation}
or, equivalently,
\begin{equation}
\d\tau_2 = {\ts\frac16}\,\bigl(|\tau_2|^2\,\sigma + \stars(\tau_2\w\tau_2)\bigr)
\end{equation}
must hold everywhere on~$M$.
\end{corollary}

\begin{proof}
Using~\eqref{eq: dtau as Ric when closed}, one obtains, after
using~\eqref{eq: G2 quartic}, the orthogonality of~$\Omega^3_1(M,\sigma)$
and~$\Omega^3_{27}(M,\sigma)$, the identity~\eqref{eq: 27-piece quartic id},
and the Cauchy-Schwartz
inequality,
\begin{equation}\label{eq: dtau3 ineq}
\begin{aligned}
\d\bigl({\ts\frac13}{\tau_2}^3\bigr)
&= {\tau_2}^2 \w \d \tau_2 
 = {\tau_2}^2 \w\bigl({\ts\frac3{14}}\,|\tau_2|^2\,\sigma
+{\ts\frac12}\stars(\tau_2\w\tau_2)
  -{\ts\frac12}\,\is\bigl(\Ric^0(g_\sigma)\bigr)\bigr)\\
& = {\ts\frac27} |\tau_2|^4\stars 1
   -{\ts\frac12}\,{\tau_2}^2\w \is\bigl(\Ric^0(g_\sigma)\bigr)\\
& = {\ts\frac27} |\tau_2|^4\stars 1
   -{\ts\frac12}\,\bigl({\tau_2}^2+{\ts\frac17}\,|\tau_2|^2\stars\sigma\bigr)
      \w \is\bigl(\Ric^0(g_\sigma)\bigr)\\
& \ge \left({\ts\frac27} |\tau_2|^4
   -{\ts\frac12}\sqrt{{\ts\frac67}}
   \,|\tau_2|^2\,\bigl|\is\bigl(\Ric^0(g_\sigma)\bigr|\right)\,\stars1.
\end{aligned}
\end{equation}
Now, the expression at the end of~\eqref{eq: dtau3 ineq} 
will be a nonnegative multiple of
the volume form~$\stars1$ as long as
\begin{equation}\label{eq: Ric0 bdd by Scal}
\bigl|\Ric^0(g_\sigma)\bigr|
=\sqrt{{{\ts\frac18}}}\,\bigl|\is\bigl(\Ric^0(g_\sigma)\bigr|
\le \sqrt{{\ts\frac1{21}}}|\tau_2|^2 = -\sqrt{{\ts\frac4{21}}}\Scal(g_\sigma).
\end{equation}

Since~$-\Scal(g_\sigma)\ge0$, the inequality~\eqref{eq: Ricci pinch ineq}
with~$C<1$ will evidently imply that the expression at the end 
of~\eqref{eq: dtau3 ineq} is a positive multiple of~$|\tau_2|^2\stars1$.
By Stokes' theorem, this will imply that~$\tau_2$ vanishes identically,
as desired.

Suppose now that \eqref{eq: Ricci pinch ineq} holds with~$C=1$.  Then
the expression at the end of~\eqref{eq: dtau3 ineq} is still a nonnegative 
multiple of~$|\tau_2|^2\stars1$ and hence, by Stokes' theorem, must
vanish identically. However, by the strong form of the Cauchy-Schwartz
inequality, this can only happen if the relation
\begin{equation}\label{eq: Ric0 as tautau}
\is\bigl(\Ric^0(g_\sigma)\bigr) 
= {\ts\frac23}\bigl(\stars(\tau_2\w\tau_2)+{\ts\frac17}|\tau_2|^2\,\sigma\bigr)
\end{equation}
holds identically on the open set where~$|\tau_2|>0$.  Now,
if the locus~$|\tau_2|=0$ has any interior, then $\Ric(g_\sigma)$
vanishes on this interior since~$\sigma$ is both closed and 
coclosed there.  Thus,~\eqref{eq: Ric0 as tautau}
holds on both the open set where~$|\tau_2|>0$ and the interior
of the locus where~$|\tau_2|=0$.  Consequently, it
must hold on all of~$M$, as desired.
\end{proof}

\begin{remark}[Extremally Ricci-pinched closed $\Gtwo$-structures]
Note that another way of phrasing Corollary~\ref{cor: pinched Ricci}
is to use \eqref{eq: dtau3 ineq} to show that the inequality
\begin{equation}\label{eq: pinched Ricci integrated}
\int_M\bigl|\Ric^0(g_\sigma)\bigr|^2\,\stars1
\ge \frac{4}{21}\int_M\Scal(g_\sigma)^2\,\stars1
\end{equation}
holds for any closed $\Gtwo$-structure~$\sigma$ 
on a compact manifold~$M^7$ and that equality holds 
in~\eqref{eq: pinched Ricci integrated} 
if and only if~$\sigma$ satisfies
\begin{equation}\label{eq: pinched Ricci defn}
\d\sigma = 0, \qquad 
\d\stars\sigma = \tau\w\sigma,\qquad
\d\tau = {\ts\frac16}\,\bigl(|\tau|^2\,\sigma + \stars(\tau\w\tau)\bigr).
\end{equation}

Indeed, Corollary~\ref{cor: pinched Ricci} suggests that the
$\Gtwo$-structures~$\sigma$ that satisfy~\eqref{eq: pinched Ricci defn}
might be of particular interest, since these are, in some sense, the most
`extremally Ricci-pinched' that a closed $\Gtwo$-structure can 
be on a compact $7$-manifold.  

One can see that there are some rather subtle restrictions on
such structures on compact manifolds by developing these equations 
a bit further:  Note that \eqref{eq: pinched Ricci defn} implies
\begin{equation}\label{eq: pinched Ricci d tau cubed}
\begin{aligned}
\d\bigl({\tau}^3\bigr)
&= 3\,{\tau}^2 \w \d \tau 
 = {\tau}^2 \w\bigl({\ts\frac12}\,|\tau|^2\,\sigma
+{\ts\frac12}\stars(\tau\w\tau)\bigr)\\
&= -{\ts\frac12}\,|\tau|^4\stars 1 
   +{\ts\frac12} |\tau\w\tau|^2\,\stars 1
= 0.
\end{aligned}
\end{equation}
On the other hand, computation using the structure equations
and \eqref{eq: pinched Ricci defn} yields
\begin{equation}\label{eq: pinched Ricci ddtau decomp}
0 = \d\bigl(\d\tau\bigr) 
=\d\left({\ts\frac16}\,\bigl(|\tau|^2\,\sigma+\stars(\tau\w\tau)\bigr)\right)
= \alpha \w \sigma + \stars\gamma
\end{equation}
where~$\gamma$ lies in~$\Omega^3_{27}(M,\sigma)$ and
\begin{equation}
\alpha 
= {\ts\frac18}\,\bigl(\d\bigl(|\tau|^2\bigr)-{\ts\frac29}\stars(\tau^3)\bigr).
\end{equation}
Consequently, any solution of~\eqref{eq: pinched Ricci defn} must satisfy%
\footnote{The vanishing of~$\gamma$ as defined in 
\eqref{eq: pinched Ricci ddtau decomp} imposes $27$ more equations
on the covariant derivative of~$\tau$, but these are not as easily
stated as~\eqref{eq: pinched Ricci d norm tau}.}
\begin{equation}\label{eq: pinched Ricci d norm tau}
\d\bigl(|\tau|^2\bigr) = {\ts\frac29}\stars(\tau^3).
\end{equation}
Combining this with~\eqref{eq: pinched Ricci d tau cubed}
yields
\begin{equation}\label{eq: norm tau harmonic}
\Delta_\sigma (|\tau|^2\bigr) = 0.
\end{equation}

Assume now that~$M$ is compact and connected.  It then follows
from~\eqref{eq: norm tau harmonic} that~$|\tau|^2$ 
must be a constant.  

Of course, if~$|\tau|^2=0$, then~$\tau=0$
and~$\sigma$ is coclosed and hence $g_\sigma$-parallel.
Thus, assume from now on that~$|\tau|^2 > 0$.

Then~\eqref{eq: pinched Ricci d norm tau} implies that~$\tau^3=0$. 
However,~$|\tau\w\tau|^2=|\tau|^4\not=0$, implying 
that~$\tau$ has constant rank~$4$
(instead of the \emph{a priori} maximum of~$6$) and hence that
$\tau\w\tau$ is a nonzero simple $4$-form of constant norm.  
Using \eqref{eq: G2 cubic express} and the fact that~$\tau^3=0$
then yields
\begin{equation}\label{eq: pinched Ricci d tau squared}
\begin{aligned}
\d\bigl({\tau}^2\bigr)
&= 2\,\tau \w \d \tau 
 = {\ts\frac13}\tau\w\bigl(\,|\tau|^2\,\sigma
        +\stars(\tau\w\tau)\bigr)\\
&= -{\ts\frac13}\,|\tau|^2\stars\tau 
   +{\ts\frac13}\,|\tau|^2\stars\tau 
= 0,
\end{aligned}
\end{equation}
So that the simple~$4$-form $\tau\w\tau$ is closed.

Since~$\tau\w\tau$ is simple with constant norm, 
the $3$-form~$\stars(\tau\w\tau)$ is also nonzero
and simple, with constant norm.  Moreover, in view of the constancy of~$|\tau|^2$, expanding $\d\bigl(\d\tau\bigr) = 0$ and using
\eqref{eq: pinched Ricci defn} shows that~$\stars(\tau\w\tau)$ 
is also closed.

Consequently, the tangent bundle of~$M$ splits
as an orthogonal direct sum of two integrable subbundles
\begin{equation}\label{eq: TM splits 3 4}
TM = P\oplus Q
\end{equation}
with~$P = \{\,v\in TM\ \vrule\ v\lhk(\tau\w\tau) =0\,\}$ of
rank~$3$ and~$Q = \{\,v\in TM\ \vrule\ v\lhk\stars(\tau\w\tau) =0\,\}$
of rank~$4$.  The $P$-leaves are calibrated 
by~$-|\tau|^{-2}\,\stars(\tau\w\tau)$ while the~$Q$-leaves are
calibrated by~$-|\tau|^{-2}\,(\tau\w\tau)$. (The reason for the minus signs
is that they correctly orient the $P$-leaves as associative
submanifolds and the~$Q$-leaves as coassociative submanifolds.)

The Ricci curvature in this case simplifies to
\begin{equation}\label{eq: Ricci in ext pinched compact case}
\Ric(g_\sigma) = {\ts\frac1{12}}\,\js\bigl(\stars(\tau\w\tau)\bigr)
= -{\ts\frac16}\,|\tau|^2\,\bigl(g_\sigma\bigr)\vrule_{P}\,,
\end{equation}
so that, in particular, the Ricci curvature is nonpositive, 
with one eigenvalue~$-{\ts\frac16}\,|\tau|^2$ of multiplicity~$3$
and the other eigenvalue~$0$ of multiplicity~$4$.

\begin{example}[A homogeneous example]\label{ex: homog ex pinched Ricci}
Just how general the $\Gtwo$-structures~$\sigma$ satisfying 
\eqref{eq: pinched Ricci defn} with~$\tau\not=0$ are, even locally, 
is an interesting question.  
I will now show that these equations do have a nontrivial solution,
by producing a (homogeneous) example.  

Let~$G$ be the group of volume-preserving
affine transformations of~$\bbC^2$.  Thus~$G$ can be regarded as
the matrix group consisting of the $3$-by-$3$  matrices with
complex entries of the form
\begin{equation}
g = \begin{pmatrix}a & b & x\\ c & d & y \\ 0 & 0 & 1 \end{pmatrix}
\end{equation}
where~$ad-bc=1$.  Write the canonical left-invariant form on~$G$ as
\begin{equation}\label{eq: left inv form of G}
\alpha = g^{-1}\,\d g = 
\begin{pmatrix}
-\omega^1 + \iC\,\eta^1 & -\omega^3-\eta^3 +\iC\,(\eta^2-\omega^2) & 
\omega^4 + \iC\,\omega^5\\
-\omega^3+\eta^3 +\iC\,(\eta^2+\omega^2) &  \omega^1 - \iC\,\eta^1 &  
\omega^6 - \iC\,\omega^7\\
0&0&0
\end{pmatrix}.
\end{equation}
Then $\d\alpha = -\alpha\w\alpha$ implies that
the left-invariant $3$-form~$\tilde\sigma$ defined by
\begin{equation}\label{eq: tilde sigma def}
\tilde\sigma = \omega^{123}+\omega^{145}+\omega^{167}
        +\omega^{246}-\omega^{257}-\omega^{347}-\omega^{356}
\end{equation}
(where~$\omega^{ijk}$ stands for the wedge 
product~$\omega^i\w\omega^j\w\omega^k$, etc.) satisfies~$\d\tilde\sigma=0$.
Consequently,~$\tilde\sigma$ is the pullback 
to~$G$ of a definite $3$-form~$\sigma$ on the left coset 
space~$M^7 = G/\SU(2)$.  (Here, ~$\SU(2)\subset G$ is the subgroup
whose left cosets are the integral leaves of the differentially 
closed system~$\omega^i=0$ on~$G$.)  Moreover, letting~$\pi:G\to M$
denote the coset projection, one sees that
\begin{equation}
\pi^*\bigl(\stars\sigma\bigr)
= \omega^{4567}+\omega^{2367}+\omega^{2345}
+\omega^{1357}-\omega^{1346}-\omega^{1256}-\omega^{1247}
\end{equation}
while
\begin{equation}
\pi^*\bigl(g_\sigma\bigr)
= (\omega^1)^2+\cdots+ (\omega^7)^2.
\end{equation}
Finally, one finds that there exists a $2$-form~$\tau$ on~$M$
so that
\begin{equation}
\pi^*(\tau) = 6\,\omega^{45} - 6\,\omega^{67}.
\end{equation}
The equation~$\d\alpha = -\alpha\w\alpha$ then implies 
that the pair~$(\sigma,\tau)$ satisfy~\eqref{eq: pinched Ricci defn}.

Note that~$M$ is diffeomorphic to~$\bbR^7$ and that the~$P$-leaves 
and~$Q$-leaves are, respectively, the fibers of 
maps~$M\to\bbC^2=G/\SL(2,\bbC)$ and~$M\to\SL(2,\bbC)/\SU(2)$.
\end{example}

\end{remark}

\begin{remark}[Natural equations for closed \Gtwo-structures]
\label{rem: natrl eqs for closed G2 strcts}
Let~$\lambda$ be a constant and consider the system of
equations
\begin{equation}\label{eq: general natrl closed G2}
\d\sigma = 0, \qquad 
\d\stars\sigma = \tau\w\sigma,\qquad
\d\tau = {\ts\frac17}\,|\tau|^2\,\sigma 
+ \lambda\,\bigl( {\ts\frac17}\,|\tau|^2\,\sigma + \stars(\tau\w\tau)\bigr).
\end{equation}
for a \Gtwo-structure~$\sigma$ on a manifold~$M^7$.
This family includes both the Einstein condition~($\lambda=\frac12$)
and the `extremally pinched Ricci' condition~($\lambda=\frac16$).
Indeed, in view of~\eqref{eq: 1 part dtau2} 
and~\eqref{eq: beta2 decomp} and since
$\Ss^2(\Vs_{0,1})\simeq\Vs_{0,0}{\oplus}\Vs_{2,0}{\oplus}\Vs_{0,2}$
while~$\L^3(\Vs_{1,0})
\simeq\Vs_{0,0}{\oplus}\Vs_{0,1}{\oplus}\Vs_{2,0}$,
the $1$-parameter family of natural 
equations~\eqref{eq: general natrl closed G2}
for closed \Gtwo-structures describes the most general way 
in which~$\d\tau$ can be prescribed naturally and quadratically 
in terms of~$\tau$.  In view of the fact that~$\d\tau$ can
have no component in~$\Omega^3_7(M,\sigma)$ and that the
component of~$\d\tau$ in~$\Omega^3_1(M,\sigma)$ is determined
by~\eqref{eq: 1 part dtau2}, it follows 
that~\eqref{eq: general natrl closed G2} is a system 
of~$27$ ($= \dim\Vs_{2,0}$) equations 
for a closed \Gtwo-structure~$\sigma$.  In view of the
discussion in Remark~\ref{rem: diff inv closed structures}, one
should regard~\eqref{eq: general natrl closed G2} as an 
overdetermined system of~\textsc{pde}.  This
system is not involutive for any value of~$\lambda$, 
as the following discussion will show.

First, the computation~\eqref{eq: pinched Ricci d tau cubed} 
can be redone for \Gtwo-structures
satisfying~\eqref{eq: general natrl closed G2}, yielding
\begin{equation}\label{eq: gen natl d tau cubed}
\d\bigl({\tau}^3\bigr) 
= \frac{3(6\lambda{-}1)}{7}\,|\tau|^4\stars 1\,.
\end{equation}
In particular, on a compact $7$-manifold,  
the only value of~$\lambda$ that is possible for such a structure 
with~$\tau$ not identically zero is~$\lambda=\frac16$.

Redoing the computation~\eqref{eq: pinched Ricci ddtau decomp} 
using the structure equations and~\eqref{eq: general natrl closed G2} 
instead of \eqref{eq: pinched Ricci defn} yields
\begin{equation}\label{eq: gen natrl ddtau decomp}
0 = \d\bigl(\d\tau\bigr) 
= \alpha \w \sigma + \stars\gamma
\end{equation}
where~$\gamma$ lies in~$\Omega^3_{27}(M,\sigma)$ and
\begin{equation}\label{eq: gen natrl ddtau decomp alpha}
\alpha 
= \frac{\lambda(2\lambda{-}1)}{4}\stars(\tau^3)
- \frac{(3\lambda{-}4)}{28}\,\d\bigl(|\tau|^2\bigr).
\end{equation}
Consequently, any solution of~\eqref{eq: general natrl closed G2} 
satisfies
\begin{equation}\label{eq: gen natrl d norm tau}
(3\lambda{-}4)\,\d\bigl(|\tau|^2\bigr) 
= 7\lambda(2\lambda{-}1)\,\stars(\tau^3).
\end{equation}

When~$\lambda=\frac43$, this condition 
implies~$\tau^3=0$, which, by~\eqref{eq: gen natl d tau cubed},  
then implies~$|\tau|=0$, i.e.,~$\tau=0$.
Thus, there are no \Gtwo-structures~$\sigma$ satisfying 
\eqref{eq: general natrl closed G2} with~$\lambda=\frac43$
except those that are closed and coclosed.

When~$\lambda\not=\frac43$,
the system~\eqref{eq: gen natrl d norm tau} represents $7$
`new' second order equations on~$\sigma$ that are not 
algebraic consequences of~\eqref{eq: general natrl closed G2}.
The existence of these `new' equations
implies that the system~\eqref{eq: general natrl closed G2}
is not involutive.  

Even beyond this, when~$\lambda\not=0$, the vanishing 
of the term~$\gamma$ in~\eqref{eq: gen natrl ddtau decomp} yields
$27$ more equations of second order on~$\sigma$ that are not
algebraic consequences of~\eqref{eq: general natrl closed G2}
and~\eqref{eq: gen natrl d norm tau} combined.  Whether further
differentiation of these combined equations would yield more 
second (or even first) order equations remains to be seen.
It is this phenomenon that makes the analysis of systems of
type~\eqref{eq: general natrl closed G2} troublesome.
\end{remark}

\section{The Torsion-free Case}
\label{sec: torfreecase} 

A \Gtwo-structure~$\sigma\in\Omega^3_+(M)$ is said to be
\emph{torsion-free} if all of its four torsion forms vanish.
There is an aspect of the geometry of torsion-free 
\Gtwo-structures that is analogous to the K\"ahler identities 
in complex Riemannian geometry and that is the concern of this section.

The material in this section was the result of a joint 
project with F. Reese Harvey and was carried out between
1991 and 1994.

\subsection{Reference modules}
It will be convenient to chose a `reference' representation for
each of the irreducible \Gtwo-modules that appear in the
exterior algebra on~$\Vs_{1,0}$.  

Given any \Gtwo-structure $\sigma\in\Omega^3_+(M)$, these
will be chosen to correspond to the
spaces of differential forms $\Omega^0(M)$, $\Omega^1(M)$,
$\Omega^2_{14}(M,\sigma)$, and $\Omega^3_{27}(M,\sigma)$.  For
simplicity, these spaces will be referred to as $\Omega_1$,
$\Omega_7$, $\Omega_{14}$, and $\Omega_{27}$ when $M$ and 
$\sigma$ are clear from context.  

\subsection{Exterior derivative identities}
When a \Gtwo-structure~$\sigma$ has vanishing intrinsic torsion,
the fundamental forms~$\sigma$ and~$\stars\sigma$ are parallel 
with respect to the natural connection (which is torsion-free) 
and so are all of the various natural isomorphisms between the
different constituents of the bundle of exterior differential
forms.  Consequently, the various differential operators that
one can define by decomposing the exterior derivative into 
its constituent components are really manifestations of first
order differential operators between the abstract bundles.
Thus, there will be identities (analogous to the identities
one proves in K\"ahler geometry) between these different 
manifestations.  In this subsection, these will be made 
explicit.  Essentially, the proof of the following proposition
is a matter of checking constants.

\begin{proposition}[Exterior derivative identities]\label{prop: G2 d identities}
Suppose that $\sigma$ is a torsion-free \Gtwo-structure 
on $M$.  Then, for all $p,q\in\{1,7,14,27\}$, 
there exists a first order differential operator 
$\d^p_q\colon\Omega_p\to\Omega_q$, so that the exterior derivative 
formulas given in Table~\ref{tab: d formulae} hold for all $f\in\Omega_1$, 
$\alpha\in\Omega_7$, $\beta\in\Omega_{14}$, and 
$\gamma\in\Omega_{27}$.
These operators are non-zero except for 
$\d^1_{27}$, $\d^{27}_1$, $\d^1_{14}$, $\d^{14}_1$, $\d^1_1$, and 
$\d^{14}_{14}$.  With respect to the natural metrics on the underlying 
bundles, $(\d^p_q)^* = \d^q_p$. 
The identity $\d^2=0$ is equivalent to the 
second order identities on the operators $\d^p_q$ listed
in Table~\ref{tab: d^2=0 formulae}. Finally, the 
formulas for the Hodge Laplacians in terms of the operators $\d^p_q$ 
are as given in Table~\ref{tab: Lapl formulae}.
\end{proposition}

\begin{proof} The operators $\d^p_q$ are defined by decomposing
the exterior derivative operator into types (much as $\partial$
and~$\bar\partial$ are defined in K\"ahler geometry by the projection
of the exterior derivative into types).  
For example, take the formula $\d^7_7\alpha = 
\stars(\d(\alpha\w\stars\sigma))$ as the \emph{definition} 
of~$\d^7_7\colon\Omega_7\to\Omega_7$ and define 
$\d^7_{27}\alpha$ to be the $\Omega^3_{27}(M,\sigma)$-component of 
$\d\bigl(\stars(\alpha\w\stars\sigma)\bigr)$.  Verifying the exterior 
derivative formulas is a routine matter that is best left to the reader.  
Once these have been established, the second order identities and the 
Laplacian formulas follow by routine computation. 
\end{proof}

\begin{table} 
$$\vcenter{\openup 1\jot \def\d#1#2{\mathrm{d}^#1_#2}
\halign{$\displaystyle #$\hfil&&\ $#$\ &\hfil$\displaystyle #$\cr
\mathrm{d}\,f  &=
&&&\d17f\cr
\mathrm{d}\,(f\,\sigma)  &=
&&&\d17f\w\sigma\cr
\mathrm{d}\,(f\,\stars\sigma)  &=
&&&\d17f\w\stars\sigma\cr
\noalign{\vskip 6pt}
\mathrm{d}\,\alpha  &=
&&&{\ts{\frac13}}\stars(\d77\alpha\w\stars\sigma)
&&+\d7{{14}}\alpha\cr
\mathrm{d}\,\stars(\alpha\w\stars\sigma) &=
&-{\ts{\frac37}}\,\d71\alpha\,\sigma
&&-{\ts{\frac12}}\stars(\d77\alpha\w\sigma)
&&
&&+\d7{{27}}\alpha\cr
\mathrm{d}\,\stars(\alpha\w\sigma) &=
&{\ts{\frac47}}\,\d71\alpha\,\stars\sigma
&&+{\ts{\frac12}}\,\d77\alpha\w\sigma
&&
&&+\stars\d7{{27}}\alpha\cr
\mathrm{d}\,(\alpha\w\sigma) &=
&&&{\ts{\frac23}}\,\d77\alpha\w\stars\sigma
&&-\stars\d7{{14}}\alpha\cr
\mathrm{d}\,(\alpha\w\stars\sigma) &=
&&&-\stars\d77\alpha\cr
\mathrm{d}\,(\stars\alpha) &=
&- \d71\alpha \,\stars1\cr
\noalign{\vskip 6pt}
\mathrm{d}\,\beta &=
&&&{\ts{\frac14}}\,\stars(\d{{14}}7\beta\w\sigma)
&&
&&+\d{{14}}{{27}}\beta\cr
\mathrm{d}\,(\stars\beta) &=
& & &\stars\d{{14}}7\beta\cr
\noalign{\vskip 6pt}
\mathrm{d}\,\gamma &=
& & &{\ts{\frac14}}\,\d{{27}}7\gamma\w\sigma
& &
&&+\stars\d{{27}}{{27}}\gamma\cr
\mathrm{d}\,(\stars\gamma) &=
& &&-{\ts{\frac13}}\,\d{{27}}7\gamma\w\stars\sigma
&&-\stars\d{{27}}{{14}}\gamma\cr
\noalign{\vskip 10pt}
} }$$
\caption{ Exterior derivative formulae} \label{tab: d formulae}
\end{table}

\begin{table}
$$\halign{&$#$\qquad\cr
&\d^7_7\,d^1_7=0\hfil&d^7_{14}\,\d^1_7=0\hfil\cr
\noalign{\vskip 10pt}
\d^7_1\,\d^7_7=0 
&\begin{aligned}\d^{14}_7\,\d^7_{14} &= {\ts{\frac23}}(\d^7_7)^2\cr
\d^{27}_7\,\d^7_{27} &= (\d^7_7)^2+{\ts{\frac{12}7}}\d^1_7\,\d^7_1\cr
\end{aligned}
&\d^7_{14}\,\d^7_7 + 2\,\d^{27}_{14}\,\d^7_{27}=0
&\begin{aligned}3\,\d^{14}_{27}\,\d^7_{14} + \d^7_{27}\,\d^7_7&=0\cr
2\,\d^{27}_{27}\,\d^7_{27} - \d^7_{27}\,\d^7_7&=0\cr\end{aligned}
\cr
\noalign{\vskip 10pt}
\d^7_1\,\d^{14}_7=0 
&\begin{aligned}\d^7_7\,\d^{14}_7 + 2\,\d^{27}_{7}\,\d^{14}_{27}&=0\cr
\end{aligned}
&{}
&\begin{aligned}
\d^{7}_{27}\,\d^{14}_{7} + 4\,\d^{27}_{27}\,\d^{14}_{27}&=0\cr
\end{aligned}
\cr
\noalign{\vskip 10pt}
&\begin{aligned}3\,\d^{14}_{7}\,\d^{27}_{14} + \d^7_{7}\,\d^{27}_7&=0\cr
2\,\d^{27}_{7}\,\d^{27}_{27} - \d^7_{7}\,\d^{27}_7&=0\cr\end{aligned}
&\begin{aligned}
\d^{7}_{14}\,\d^{27}_{7} + 4\,\d^{27}_{14}\,\d^{27}_{27}&=0\cr
\end{aligned}
&{}\cr
\noalign{\vskip 10pt}
}$$
\caption{Second order identities}\label{tab: d^2=0 formulae}
\end{table}

\begin{table}
\begin{equation*}
\begin{aligned}
\Delta\,f 
 &= \d^7_1\,\d^1_7\,f\cr
\Delta\,\alpha 
&= \left((\d^7_7)^2 + \d^1_7\,\d^7_1\right)\,\alpha\cr
\Delta\,\beta 
&= \left({\ts{\frac54}}\d^7_{14}\,\d^{14}_7 + 
\d^{27}_{14}\,\d^{14}_{27}\right)\,\beta\cr
\Delta\,\gamma 
&= \left({\ts{\frac7{12}}}\d^7_{27}\,\d^{27}_7 + 
  \d^{14}_{27}\,\d^{27}_{14}+(\d^{27}_{27})^2\right)\,\gamma\cr
\noalign{\vskip 10pt}
\end{aligned}
\end{equation*}
\caption{Laplacians}\label{tab: Lapl formulae}
\end{table}

\begin{remark}[Torsion perturbations]\label{rem: torsion perturbations}
In the general case of a \Gtwo-structure
with torsion, all of the formulae in the tables listed above must
be modified by lower order terms.  For example, 
in Table~\ref{tab: d formulae} the second line would be modified
to
\begin{equation}
\d(f\,\sigma) = \d^1_7f\w\sigma 
+ f\tau_0\,\stars\sigma + 3f\tau_1\w\sigma + f\,\stars\tau_3.
\end{equation}
The zero right hand sides in Table~\ref{tab: d^2=0 formulae} have
to be replaced by first order operators whose coefficients depend
on the torsion terms and, in Table~\ref{tab: Lapl formulae}, one
must take into account which particular part of the exterior
algebra a given form occupies before writing down the appropriate
formula for the Laplacian.  It is not true, in general, that 
$\Delta (f\,\sigma) = \Delta f\,\sigma$, for example. 
\end{remark}

\section{Deformation and Evolution of \Gtwo-structures} 
\label{sec: defandevol}

The material in this section was the result of a joint
project with Steve Altschuler and was carried out between
1992 and 1994.  Our goal was to understand the long time
behavior of the Laplacian heat flow defined below for
closed \Gtwo-structures on compact $7$-manifolds, specifically,
to understand conditions under which one could prove that
this flow converged to a \Gtwo-structure that is both
closed and coclosed.  Nowadays, this flow is called
the \emph{Hitchin flow} after Hitchin's fundamental 
paper~\cite{MR02m:53070}.

We were never able to prove long-time existence
under any reasonable hypotheses, so we wound up not
publishing anything on the subject, although we did get some
interesting results and formulae that I have not seen
so far in the literature.%
\footnote{I would be happy to 
learn of any places where these results have appeared
so that I can properly acknowledge them in future versions
of this article.}

\subsection{The deformation forms}\label{ssec: deform forms}
It turns out to be quite easy to describe deformations
of \Gtwo-structures.  The following result is well-known
and can be found most explicitly in Joyce's 
treatment~\cite[\S10.3]{MR01k:53093}, though the notation
is somewhat different.  It is included here to establish
notation for the discussion to follow.

\begin{proposition}[Deformation forms]\label{prop: deformation forms}
Let~$\sigma_t\in\Omega^3_+(M)$ be a smooth
$1$-parameter family of $\Gtwo$-structures on~$M$.  Let~$g_t$
and~$\ast_t$ denote the underlying metric and Hodge star
operator associated to~$\sigma_t$.  Then there
exist three differential forms~$f^0_t\in\Omega^0(M)$, 
$f^1_t\in\Omega^1(M)$, 
and~$f^3_t\in\Omega^3_{27}(M,\sigma_t)\subset\Omega^3(M)$ that
depend differentiably on~$t$ and that are uniquely
characterized by the equation {\upshape(}in which the $t$-dependence
has been suppressed for notational clarity{\upshape)}
\begin{equation}\label{eq: f deform defined}
\frac{\d\hfil}{\d t}\bigl(\sigma\bigr)
= 3 f^0\,\sigma + \stars\bigl(f^1\w\sigma\bigr) + f^3\,.
\end{equation}
Moreover, the associated metric and dual $4$-forms satisfy
\begin{equation}\label{eq: metric deform formula}
\frac{\d\hfil}{\d t}\bigl(g\bigr)
= 2 f^0\,g + {\ts\frac12}\,\js\bigl(f^3\bigr)
\end{equation}
and 
\begin{equation}\label{eq: dual 4form deform formula}
\frac{\d\hfil}{\d t}\bigl(\stars\sigma\bigr)
= 4 f^0\,\stars\sigma + f^1\w\sigma -\stars f^3\,.
\end{equation}
\end{proposition}

\begin{definition}[The deformation forms]\label{def: def forms}
The forms~$f^0_t$, $f^1_t$, and~$f^3_t$ associated to the family~$\sigma_t$
will be referred to as the \emph{deformation forms} of the family.
\end{definition}

One immediate consequence of Proposition~\ref{prop: deformation forms}
is a formula for the variation of the volume form:
\begin{equation}\label{eq: vol form deform formula}
\frac{\d\hfil}{\d t}\bigl(\stars1\bigr)
= 7 f^0\,\stars1\,.
\end{equation}

By the same techniques, one can derive a second order expansion:

\begin{proposition}[Taylor expansion formula]\label{prop: Taylor2} 
Let~$\phi\in\Omega^3_+(M)$ be a $\Gtwo$-structure.  Then for all
$b_0\in\Omega^0(M)$, $b_1\in\Omega^1(M)$, and~$b_3\in\Omega^3_{27}(M,\phi)$
of sufficiently small $C^0$-norm, the $3$-form
\begin{equation}\label{eq: linear deform}
\sigma = \phi + \left(3b_0\,\phi + \starf(b_1\w\phi) + b_3\right)
\end{equation}
is definite.  Moreover, there is an expansion of the form 
\begin{equation}
\begin{aligned}
\stars\sigma 
&= \starf\phi + \bigl(4b_0\,\starf\phi + b_1\w\phi - \starf b_3\bigr) 
+ \left(\,2(b_0)^2 + {\ts\frac2{21}}\,|b_1|^2_\phi - {\ts\frac1{42}}\,|b_3|^2_\phi\,\right) \starf\phi \\
&\qquad\qquad + Q_1(b_0,b_1,b_3)\w\phi + \starf Q_3(b_0,b_1,b_3) 
+ R(b_0,b_1,b_3)
\end{aligned}
\end{equation}
where $Q_1$~{\upshape(}a $1$-form{\upshape)} 
and $Q_3$~{\upshape(}a $3$-form in~$\Omega^3_{27}(M,\phi)${\upshape)} 
are quadratic in the coefficients of the~$b_i$ 
and $R$ is a $4$-form that vanishes to order~$3$ 
in the coefficients of the~$b_i\,$.
Consequently, there is an expansion of the form
\begin{equation}\label{eq: volume expansion}
\stars1 = \left(1 +  7b_0 
      + \bigl(\,14(b_0)^2 + {\ts\frac23}|b_1|^2_\phi 
    -{\ts\frac16}|b_3|^2_\phi\,\bigr)
        + r(b_0,b_1,b_3)\,\right)\,\starf1 
\end{equation}
where $r$ vanishes to order~$3$ in~$(b_0,b_1,b_3)$.
\end{proposition}

\subsection{The Laplacian evolution}\label{ssec: Lap evol}
A natural evolution equation for~\Gtwo-structures is
the (nonlinear) Laplacian evolution equation 
for~$\sigma\in\Omega^3_+(M)$ defined as follows:
\begin{equation}\label{eq: Lap evol}
\frac{d\hfil}{dt}(\sigma) = \Delta_\sigma \sigma\,.
\end{equation}
This equation is diffeomorphism invariant and hence cannot be
elliptic in the strict sense.  However, it is not difficult to
compute the linearization and see that it is transversely
elliptic, i.e., elliptic transverse to the action of the 
diffeomorphism group.

Thus, the by-now standard methods of DeTurck and Hamilton
can be applied to show that, if~$M$ is compact, then for any smooth
$\sigma_0\in\Omega^3_+(M)$ there exists an extended 
number~$T$ satisfying~$0<T\le\infty$ 
and a $1$-parameter family~$\sigma(t)\in\Omega^3_+(M)$ defined 
for all~$t$ such that~$0<t<T$ so that the family 
satisfies~\eqref{eq: Lap evol} and so that~$\sigma(t)$ 
approaches~$\sigma_0$ uniformly as~$t$ approaches~$0$ from above.
The fundamental issue then becomes to understand the behavior
of the family as~$t$ approaches~$T$. 

For general~$\sigma$, the formula for the Laplacian in terms
of the torsion forms is not too illuminating:
\begin{equation}
\Delta_\sigma \sigma = \d\bigl(\tau_2-4\,\tau_1^\sharp\lhk\sigma\bigr)
  + \stars\d\bigl(\tau_0\,\sigma+3\,\tau_1^\sharp\lhk\stars\sigma+\tau_3\bigr).
\end{equation}
This can be further expanded, but the general formula
becomes unwieldy rather quickly.

\subsubsection{Evolution of closed forms}\label{ssec: closed evol}
Suppose now that the initial form~$\sigma$ is closed, i.e., that
$\tau_0$, $\tau_1$ and~$\tau_3$ are all zero initially.  It
is not difficult to show that the Laplacian flow preserves this
condition, i.e., that the family~$\sigma(t)$ consists of closed forms.%

For notational simplicity, for the rest of this section, $\tau_2$
will be denoted simply as~$\tau$.  Also, in the calculations to follow,
$t$ will be treated as a parameter, i.e., I will regard $\d t$ as
zero when computing exterior derivatives.  Thus, the assumptions
are that
\begin{equation}\label{eq: def tau evol}
\begin{aligned}
\d\sigma & = 0\\
\d\stars\sigma & = \tau\w\sigma
\end{aligned}
\end{equation}
and that
\begin{equation}\label{eq: sigma evol}
\frac{d\hfil}{dt}(\sigma) = \d\tau.
\end{equation}
As has already been shown in~\eqref{eq: 1 part dtau2},
\begin{equation}\label{eq: dtau evol}
d\tau = {\ts\frac17}|\tau|^2\,\sigma + \gamma
\end{equation}
for some~$\gamma\in\Omega^3_{27}(M,\sigma)$.  In particular,
it follows from Proposition~\ref{prop: deformation forms} that
\begin{equation}\label{eq: starsigma evol}
\frac{d\hfil}{dt}(\stars\sigma) 
= {\ts\frac4{21}}|\tau|^2\,\stars\sigma -\stars\gamma
= {\ts\frac1{3}}|\tau|^2\,\stars\sigma -\stars\d\tau\,.
\end{equation}
Moreover, \eqref{eq: vol form deform formula} now becomes
\begin{equation}\label{eq: dvol evol}
\frac{d\hfil}{dt}(\stars1) 
= {\ts\frac13}|\tau|^2\,\stars1.
\end{equation}
In particular, note that the associated volume form~$\stars1$
is \emph{pointwise} increasing.\footnote{In view of 
Hitchin's interpretation of this flow as the gradient flow
of the volume functional on the space $[\phi]_+$, it is to be
expected that the integral of~$\stars1$ over~$M$ is increasing.}

Finally, combining~\eqref{eq: sigma evol} with the formulae
\eqref{eq: dtau as Ric when closed} and \eqref{eq: metric deform formula}, 
one gets the evolution of the metric~$g_\sigma$
in the form
\begin{equation}\label{eq: g evol}
\frac{d\hfil}{dt}(g_\sigma) 
= -2\Ric(g_\sigma) + {\ts\frac16}|\tau|^2\, g_\sigma 
                   + {\ts\frac14}\,\js\bigl(\stars(\tau\w\tau)\bigr).
\end{equation}

\begin{remark}[Hitchin's interpretation]\label{rem: HitchinInterp}
Hitchin~\cite{MR02m:53070} has given the following interpretation 
of this flow.  Suppose that~$\phi$ is a closed definite $3$-form
and on a compact $7$-manifold~$M$. Let
\begin{equation}
[\phi]_+
= \{\ \phi + \d\beta\in \Omega^3_+(M)\ \mid\ \beta\in\Omega^2(M)\ \}
\end{equation}
be the open set in the cohomology class~$[\phi]
= \{\phi+\d\beta \,\mid\,\beta\in\Omega^2(M)\}$ that consists of 
definite $3$-forms.  

Define the volume function~$V:[\phi]_+\to\bbR^+$
by~$V(\sigma) = \int_M \stars 1>0$ for~$\sigma\in[\phi]_+$.
Hitchin shows that~$\sigma\in[\phi]_+$ is a critical
point of~$V$ if and only if~$\sigma$ is coclosed (as well as closed)
and he shows that the flow~\eqref{eq: Lap evol} is the gradient
flow of the functional~$V$ (in the $L^2$ metric on~$[\phi]_+$).

Suppose that~$\phi$ is a critical point of~$V$, i.e., that $\starf\phi$
is closed.  Then by Hodge theory there is a direct sum decomposition
\begin{equation}
\d\bigl(\Omega^2(M)\bigr) 
= \{\,\Lie_Z\phi\,\mid\, Z\in\Vect(M)\,\}
  \oplus 
   \{\,\d\beta\,\mid\,\beta\in\Omega^2_{14}(M,\phi),\ \d^{14}_7\beta=0\,\}.
\end{equation}
The first summand is the tangent space to the orbit of~$\phi$
under~$\Diff^\circ(M)$ (i.e., the diffeomorphisms of~$M$ that
act trivially on~$H^*(M)$), while the second summand represents the tangent
space at~$\Diff^\circ(M){\cdot}\phi$ to the `moduli space'~$\Diff^\circ(M)\backslash[\phi]_+$.
If~$\beta\in\Omega^2_{14}$ satisfies~$d^{14}_7\beta=0$, then, 
setting~$\sigma= \phi + t\,\d\beta = \phi + t\,d^{14}_{27}\beta$, 
one finds, by~\eqref{eq: volume expansion}, that
\begin{equation}
\stars1 
= \bigl(1 -{\ts\frac16}\left|d^{14}_{27}\beta\right|^2+t^3R(t,\d\beta)\bigr)
\starf1.
\end{equation}
for some smooth remainder term~$R(t,\d\beta)$.  By the formulae
in~Table~\ref{tab: d formulae}, the equations
$\d^{14}_7\beta = \d^{14}_{27}\beta=0$
for~$\beta\in\Omega^2_{14}(M,\phi)$ imply that~$\d\beta=\delta\beta=0$.
It follows that the Hessian of~$V$ at~$\phi$ is negative definite
on~$\{\,\d\beta\,\mid\,\beta\in\Omega^2_{14}(M,\phi),\ \d^{14}_7\beta=0\,\}$.
Thus, $\Diff^\circ(M){\cdot}\phi$ is a local maximum of~$V$
on the moduli space~$\Diff^\circ(M)\backslash[\phi]_+$.%
\footnote{Hitchin
says that~$V$ is a `Morse-Bott' functional on~$[\phi]_+$, i.e.,
that~$V$ has nondegenerate critical points on the moduli space.}

In particular, it seems reasonable to expect that, for~$\sigma\in[\phi]_+$
`sufficiently near'~$\phi$ in a appropriate norm, 
the $V$-gradient flow~\eqref{eq: Lap evol} with~$\sigma$ as initial value 
would converge to a point on~$\Diff^\circ(M){\cdot}\phi$.
\end{remark}

\begin{remark}[Nonconvergence]\label{rem: poss nonconverge}
A more likely difficulty, it seems, is posed by the possibility
that there may be torsion-free \Gtwo-structures~$\phi$ for which
the volume functional is not bounded above on~$[\phi]_+$, 
so one would not expect the Laplacian flow to converge for most closed 
\Gtwo-structures in~$[\phi]_+$.

\begin{example}[Fern\'andez' closed $\Gtwo$-solvmanifold]
\label{ex: Fernandez}
Fern\'andez~\cite{88k:53070,89f:53058} has constructed compact
$7$-dimensional manifolds~$M^7$ that support a closed 
$\Gtwo$-structure~$\phi$ but that cannot, for topological reasons, 
support a torsion-free~$\Gtwo$-structure. Thus, in these
cases, the above flow cannot converge, since there will be no
critical points of~$V$ on~$[\phi]_+$.

It is instructive to look at one of her examples:  Let~$G\subset\GL(5,\bbR)$
be the subgroup that consists of matrices of the form
\begin{equation}\label{eq: g in G}
g = \begin{pmatrix}
1&0& x^2 & x^4 & x^6\\
0&1& x^3 & x^5 & x^7\\
0&0&  1  &  0  & x^1\\
0&0&  0  &  1  &  0 \\
0&0&  0  &  0  &  1\\
\end{pmatrix}
\end{equation}
where~$x^i\in\bbR$ for~$1\le i\le 7$.  Write the left-invariant 
form on~$G$ in the form
\begin{equation}\label{eq: l inv form on G}
g^{-1}\,\d g = \begin{pmatrix}
0&0& \omega^2 & \omega^4 & \omega^6\\
0&0& \omega^3 & \omega^5 & \omega^7\\
0&0&  0  &  0  & \omega^1\\
0&0&  0  &  0  &  0 \\
0&0&  0  &  0  &  0\\
\end{pmatrix}
\end{equation}
where~$\d\omega^i=0$ for~$1\le i\le 5$ while~$\d\omega^6 = \omega^1\w\omega^2$
and~$\d\omega^7 = \omega^1\w\omega^3$.

Let~$\Gamma = G\cap \GL(5,\bbZ)$ and note that~$\Gamma$ is a co-compact
discrete subgroup of~$G$.  Let~$M^7 = \Gamma\backslash G$ be the space
of right cosets of~$\Gamma$ in~$G$.  Then the~$\omega^i$ are well-defined
on~$M$ and it is easy to verify that the $3$-form
\begin{equation}
\sigma = \omega^{123}+\omega^{145}+\omega^{167} 
+\omega^{246} - \omega^{257} - \omega^{347} - \omega^{356}
\end{equation}
is a closed $\Gtwo$-structure on~$M$.  It is not coclosed, but satisfies
\begin{equation}
\d\stars\sigma = \bigl(\omega^{27}-\omega^{36}\bigr)\w\sigma.
\end{equation}
Since
\begin{equation}
\d\bigl(\omega^{27}-\omega^{36}\bigr) = 2\,\omega^{123},
\end{equation}
it follows that the flow satisfies
\begin{equation}
\sigma(t) = \E^{2t}\omega^{123}+\omega^{145}+\omega^{167} 
+\omega^{246} - \omega^{257} - \omega^{347} - \omega^{356}.
\end{equation}
The associated metric is
\begin{equation}
g(t) = \E^{4t/3}\bigl((\omega^1)^2+(\omega^2)^2+(\omega^3)^2\bigr)
+ \E^{-2t/3}\bigl((\omega^4)^2+(\omega^5)^2+(\omega^6)^2+(\omega^7)^2\bigr).
\end{equation}
In particular, note that, under this flow (which exists for all time,
both past and future), the volume of the metric increases without bound.

By the way, $M$ cannot carry a metric with holonomy a subgroup of~$\Gtwo$
for the following reason:  As Fern\'andez shows, the first Betti
number of~$M$ is~$5$.  If there were a metric~$g$ on~$M$ with holonomy
in~$\Gtwo$, then it would be Ricci-flat and hence the harmonic representatives
of the first cohomology group would give five linearly independent
$g$-parallel $1$-forms on~$M$.  However, this would imply that the
holonomy of~$M$ is trivial, which would imply that there exist seven
linearly independent parallel $1$-forms on~$M$, which would in turn imply
that the first Betti number was at least~$7$.
\end{example}

\end{remark}

\subsubsection{Further calculations}\label{sssec: further calcs}
Return now to the flow of a general closed $\Gtwo$-structure.
Taking the exterior derivative of~\eqref{eq: starsigma evol} yields
\begin{equation}\label{eq: dtausigma evol}
\frac{d\hfil}{dt}(\tau\w\sigma) 
= {\ts\frac1{3}}\d\bigl(|\tau|^2\bigr)\w\stars\sigma
  + {\ts\frac1{3}}|\tau|^2\,\tau\w\sigma -\d\stars\d\tau\,.
\end{equation}
Expanding the left hand side of this equation and using
\eqref{eq: sigma evol} yields
\begin{equation}
\frac{d\hfil}{dt}(\tau) \w \sigma 
= {\ts\frac1{3}}\d\bigl(|\tau|^2\bigr)\w\stars\sigma
  + {\ts\frac1{3}}|\tau|^2\,\tau\w\sigma -\d\stars\d\tau 
  - \tau\w\d\tau\,.
\end{equation}
(This equation can be solved for the time-derivative 
of~$\tau$ since wedging with~$\sigma$ is an isomorphism 
between~$\Omega^2$ and~$\Omega^5$.)
Recalling that~$\tau\w\stars\sigma=0$ and~$\tau\w\tau\w\sigma
= -|\tau^2|\stars1$, this yields
\begin{equation}
\frac{d\hfil}{dt}(\tau) \w \tau \w \sigma 
= - {\ts\frac1{3}}|\tau|^4\,\stars1 -\tau \w \d\stars\d\tau 
  - \tau \w \tau\w\d\tau\,.
\end{equation}
Finally, this can be used in the following computation
\begin{equation}\label{eq: d|tau|^2vol evol}
\begin{aligned}
\frac{d\hfil}{dt}(|\tau^2|\,\stars1)
&= \frac{d\hfil}{dt}(-\tau \w \tau \w \sigma) 
= -2 \frac{d\hfil}{dt}(\tau) \w \tau \w \sigma -\tau\w\tau\w\d\tau\\
&=  {\ts\frac2{3}}|\tau|^4\,\stars1 + 2 \tau \w \d\stars\d\tau
     + \tau \w \tau\w\d\tau\\
&= \bigl({\ts\frac2{3}}|\tau|^4-2|\d\tau|^2\bigr)\,\stars1
     +\d\left(2\tau\w \stars\d\tau + {\ts\frac1{3}}\tau^3\right).
\end{aligned}
\end{equation}
Integrating this equation over~$M$ yields
\begin{equation}
\frac{d\hfil}{dt}\int_M |\tau^2|\,\stars1 = 
\int_M \bigl({\ts\frac2{3}}|\tau|^4-2|\d\tau|^2\bigr)\,\stars1.
\end{equation}

This equation can be rewritten by using~\eqref{eq: dtau as Ric when closed}, 
which yields
\begin{equation}\label{eq: 2nd der expansion 1}
\frac{d\hfil}{dt}\int_M (|\tau^2|\,\stars1)
= \int_M \bigl({\ts\frac2{3}}|\tau|^4
-2\left|{\ts\frac3{14}}\,|\tau|^2\,\sigma
{+}{\ts\frac12}\stars(\tau\w\tau) 
{-}{\ts\frac12}\is\bigl(\Ric^0(g_\sigma)\bigr)\right|^2\bigr)\,\stars1.
\end{equation}

Now, going back to~\eqref{eq: dtau3 ineq} and integrating this over~$M$
yields
\begin{equation}
0 = \int_M {\ts\frac27} |\tau_2|^4\,\stars 1
   -{\ts\frac12}\,\bigl({\tau_2}^2+{\ts\frac17}\,|\tau_2|^2\stars\sigma\bigr)
      \w \is\bigl(\Ric^0(g_\sigma)\bigr),
\end{equation}
i.e.,
\begin{equation}
\int_M 
\left\langle\stars\bigl({\tau_2}^2+{\ts\frac17}\,|\tau_2|^2\stars\sigma\bigr)
      , \is\bigl(\Ric^0(g_\sigma)\bigr)\right\rangle\,\stars 1
= {\frac47}\int_M  |\tau_2|^4\,\stars 1.
\end{equation}

Using this relation and the algebraic 
identities~\eqref{eq: 27-piece quartic id}
and~\eqref{eq: beta2 decomp}, one sees that
\eqref{eq: 2nd der expansion 1} can be rewritten
in the form
\begin{equation}\label{eq: second vol der}
\frac{d\hfil}{dt}\int_M (|\tau^2|\,\stars1)
= 4\int_M \left({\ts\frac{11}{21}}\Scal(g_\sigma)^2
   -\bigl|\Ric^0(g_\sigma)\bigr|^2\,\right)\,\stars1.
\end{equation}

This equation is suggestive.  One of the reasons for wanting to study
the Laplacian flow on closed \Gtwo-structures is that it might provide
a means of constructing metrics with holonomy \Gtwo\ by starting with
a closed \Gtwo-structure~$\sigma\in\Omega^3_+(M)$ with `sufficiently
small' torsion and then running the Laplacian flow to move it closer
to a \Gtwo-structure that is both closed and coclosed.  

However, if such a procedure is to work, then the volume function
along the flow line must approach a constant and one would certainly 
expect the second derivative to become negative if the volume were
to approach the `local maximum' target volume.  However, 
\eqref{eq: second vol der} shows that, in this case, the relative
separation of the eigenvalues of the Ricci tensor 
cannot decrease too much during the flow.  This `forced separation'
is somewhat stronger than the separation implied by 
Corollary~\ref{cor: pinched Ricci}.

\bibliographystyle{hamsplain}

\providecommand{\bysame}{\leavevmode\hbox to3em{\hrulefill}\thinspace}

\end{document}